\documentclass[12pt,twoside,reqno]{amsart}
\usepackage{times}
\usepackage{amsmath,amssymb,amsthm}
\usepackage{color}
\usepackage{enumerate}
\usepackage{cite}

\allowdisplaybreaks

\newcommand{\R}{\mathbb R}
\newcommand{\N}{\mathbb N}

\newcommand{\p}{\partial}
\newcommand{\ve}{\varepsilon}
\newcommand{\f}{\frac}
\newcommand{\na}{\nabla}

\renewcommand{\th}{\theta}

\newcommand{\G}{\Gamma}

\newcommand{\ds}{\displaystyle}
\newcommand{\dbar}{{\mathchar'26\mkern-12mu d}}

\numberwithin{equation}{section}

\theoremstyle{plain}
\newtheorem{theorem}{Theorem}[section]
\newtheorem{lemma}[theorem]{Lemma}
\newtheorem{proposition}[theorem]{Proposition}

\theoremstyle{definition}

\newtheorem{corollary}[theorem]{Corollary}

\theoremstyle{remark}

\newtheorem{remark}[theorem]{Remark}

\title[minimal regularity]{Minimal regularity solutions of semilinear
  \\ generalized Tricomi equations}

\author[Z.-P.~Ruan]{Zhuoping Ruan}
\address{Department of Mathematics, Nanjing University, Nanjing
  210093, China}
\email{zhuopingruan@nju.edu.cn}

\author[I.~Witt]{Ingo Witt}
\address{Mathematical Institute, University of G\"{o}ttingen,
  Bunsenstr.~3-5, D-37073 G\"{o}ttingen, Germany}
\email{iwitt@uni-math.gwdg.de}

\author[H.-C.~Yin]{Huicheng Yin}
\address{School of Mathematical Sciences, Nanjing Normal University,
  Nanjing 210023, China}
\email{huicheng@nju.edu.cn, 05407@njnu.edu.cn}

\date{}

\keywords{Generalized Tricomi equation, minimal regularity, Fourier
  integral operators, Stricharz estimates}

\subjclass[2010]{Primary: 35L70; Secondary: 35L65, 35L67, 76N15}

\thanks{The first and third authors were supported by the NSFC
  (No.~11401299, No.~11571177) and by the Priority Academic Program
  Development of Jiangsu Higher Education Institutions.}

\begin{document}

\begin{abstract}
We prove the local existence and uniqueness of minimal regularity
solutions $u$ of the semilinear generalized Tricomi equation $\p_t^2 u
-t^m \Delta u =F(u)$ with initial data $(u(0, \cdot), \p_t u(0,
\cdot)) \in \dot{H^{\gamma}}(\R^n) \times \dot{H}^{\gamma-\frac 2
  {m+2}}(\R^n)$ under the assumption that $ |F(u) | \lesssim
|u|^\kappa$ and $|F'(u)| \lesssim |u|^{\kappa -1}$ for some
$\kappa>1$. Our results improve previous results of M.~Beals
\cite{Beals} and of~ourselves \cite{RWY1,RWY2,RWY3}.  We establish
Strichartz-type estimates for the linear generalized Tricomi operator
$\p_t^2 -t^m \Delta$ from which the semilinear results are derived.
\end{abstract}

\maketitle


\section{Introduction}

In this paper, we are concerned with the local well-posedness problem
for minimal regularity solutions $u$ of the semilinear generalized
Tricomi equation
\begin{equation} \label{1.1}
\left\{ \enspace
\begin{aligned}
& \p_t^2 u -t^m \Delta u =F(u) &&\text{in }   (0,T) \times \R^n, \\
& u (0, \cdot)=\varphi  \in \dot{H}^\gamma(\R^n), \\
& \p_t u (0, \cdot)=\psi \in \dot{H}^{\gamma-\frac 2 {m+2}}(\R^n),
\end{aligned}
\right.
\end{equation}
where $n\ge 2$, $m\in\N$, $\gamma \in \R$,
$\Delta=\sum_{i=1}^n\p_i^2$, and $T>0$. The nonlinearity $F \in
C^1(\R)$ obeys the estimates
\begin{equation}\label{1.2}
  |F(u)| \lesssim |u|^{\kappa}, \quad
  |F'(u)| \lesssim |u|^{\kappa -1}
\end{equation}
for some $\kappa>1$. For $n\geq 3$ and $\kappa > \kappa_3$ (see below)
we further assume that $\kappa\in\N$ and $F(u)=\pm u^\kappa$.

Our main objective of this paper is to find the minimal number
$\gamma$ for which Eq.~\eqref{1.1} under assumption~\eqref{1.2}
possesses a unique local solution $u \in C([0, T],
\dot{H}^\gamma(\R^n))\cap L^s((0,T);L^q(\R^n))$ for certain $s,\,q$
with $\min\{s,q\}\geq\kappa$. Then $F(u)\in
L^{s/\kappa}((0,T);L^{q/\kappa}(\R^n))\subseteq
L_{\text{loc}}^1((0,T)\times\R^n)$, and Eq.~\eqref{1.1}
holds in distributions.


\smallskip

We first introduce notation used throughout this
paper. Set
\[
\mu_\ast =\frac {(m+2)n +2}{2}, \quad \kappa_\ast =\frac {\mu_\ast
  +2}{\mu_\ast -2} =\frac{(m+2)n+6}{(m+2)n-2},
\]
\begin{align*}
\kappa_0&= 1+ \frac{6\mu_\ast+m}{\mu_\ast (m+2) n} \hspace{87pt}
\text{if $n\geq3$ or $n=2$, $m\geq3$,} \\
\kappa_1&=
\begin{cases}
  2  & \ \ \text{if $n=2$, $m=1$,} \\
  \ds \frac{ (\mu_\ast +2) (m+2)(n-1)+8}{ (\mu_\ast -2)(m+2)(n-1)+8} &
  \ \ \text{if $n\geq3$ or $n=2$, $m\geq2$,}
\end{cases} \\
\kappa_2 &=\frac{\mu_\ast (\mu_\ast +2)(n-1) -2(n+1)}{\mu_\ast
  (\mu_\ast -2)(n-1) -2(n+1)},\\
\intertext{and}
\kappa_3 &=\frac{\mu_\ast-m}{\mu_\ast -m -4} \quad \hspace{98pt}
\text{if $n \ge 3$.}
\end{align*}
Note that $\mu_\ast$ is the homogeneous dimension of the degenerate
differential operator $\p_t^2-t^{m}\Delta$ and $\kappa_\ast$ is the
power $\kappa$ for which the equation $\p_t^2 u -t^m \Delta u =
\pm\,|u|^{\kappa-1}u$ is conformally invariant.
Note further that $1<\kappa_0<\kappa_1<\kappa_* < \kappa_2 <\kappa_3$
whenever it applies.

\smallskip

Now we state the main results of this paper.

\begin{theorem} \label{thm1.1}
Let $n\geq2$ and $F$ be as above. Suppose further $\kappa>\kappa_1$
and $\left(\varphi,\psi\right)\in \dot{H}^\gamma(\mathbb R^n)\times
\dot{H}^{\gamma-\frac 2 {m+2}}(\mathbb R^n)$, where
\begin{equation} \label{1.3}
\gamma = \gamma(\kappa, m, n)=
\begin{cases}
\ds\frac {n+1} 4 -\frac{n+1}{\mu_\ast (\kappa-1)}-\frac m {2 \mu_\ast
  (m+2)} & \textup{if } \kappa_1 < \kappa \le \kappa_\ast, \\
\ds\frac n 2 -\frac 4 {(m+2)(\kappa-1)} & \textup{if } \kappa\ge
\kappa_\ast.
\end{cases}
\end{equation}
Then problem \eqref{1.1} possesses a unique solution
\[
u  \in  C([0, T];  \dot{H}^\gamma(\R^n)) \cap L^s ((0, T);L^q(\R^n))
\]
for some $T>0$, where
\begin{equation}\label{1.4}
\|u\|_{C([0, T]; \dot{H}^\gamma(\R^n))}  + \|u\|_{L^s ((0, T);L^q(\R^n))}
\lesssim  \|\varphi\|_{\dot{H}^{\gamma}(\R^n)} +
\|\psi\|_{\dot{H}^{\gamma-\frac 2 {m+2}}(\R^n)}
\end{equation}
and  $q=\mu_\ast\left(\kappa-1\right)/2$,
\begin{equation*}
\frac1 s= \begin{cases}
\ds \frac{(m+2)(n-1)}4 \left( \frac 1 2 -\frac 1 q\right) +\frac m {4 \mu_\ast} &
\textup{if } \kappa_1 < \kappa \le \kappa_\ast, \\
\ds 1/q & \textup{if } \kappa\geq \kappa_*.
\end{cases}
\end{equation*}
\end{theorem}

\smallskip

\begin{remark}\label{rem1.2}
As a byproduct of the proof of Theorem~\ref{thm1.1}, we see that
problem \eqref{1.1} admits a unique global solution $u \in C([0,
  \infty); \dot{H}^\gamma(\R^n) )\cap L^\infty((0, \infty);
  \dot{H}^\gamma(\R^n) )\cap L^{\frac{\mu_\ast(\kappa-1)}2} ( \R_+
  \times \R^n )$ in case $n\geq2$, $\kappa \ge \kappa_\ast$ if
  $(\varphi, \psi)=\ve\left( u_0, u_1\right)$, $(u_0, u_1)\in
  \dot{H}^\gamma(\R^n) \times \dot{H}^{\gamma-\frac 2 {m+2}}(\R^n)$,
  and $\ve>0$ is small (cf.~\ref{case_2} and \ref{case_3} in the proof
  of Theorem~\ref{thm1.1} below). With a
  different argument, the global result $u \in L^{\frac{\mu_\ast
      (\kappa-1)} 2}( \R_+ \times \R^n)$ for problem \eqref{1.1} was
  obtained in \cite{He1}.
\end{remark}

\begin{remark}\label{rem1.4}
For $\gamma < \frac n 2 -\frac 4{(m+2)(\kappa-1)}$, one obtains
ill-posedness for problem \eqref{1.1} by scaling. More specifically,
if $u=u(t,x)$ solves the Cauchy problem \eqref{1.1}, where
$F(u)=\pm\,|u|^{\kappa-1}u$, then
\[
u_\varepsilon(t, x)= \varepsilon^{-\frac 2{\kappa -1}} u(
\varepsilon^{-1} t, \varepsilon^{-\frac {m+2} 2} x), \quad \ve>0,
\]
also solves \eqref{1.1}, with
$u_\varepsilon(0,x)=\varphi_\varepsilon(x)$, $\p_t u_\varepsilon(0,x)=
\psi_\varepsilon(x)$ for some resulting $\varphi_\varepsilon$,
$\psi_\varepsilon$. Observe that
\[
\frac {\|\varphi_\varepsilon\|_{\dot{H}^\gamma(\R^n)}}{
  \|\varphi\|_{\dot{H}^\gamma(\R^n)}} =\frac
      {\|\psi_\varepsilon\|_{\dot{H}^\gamma(\R^n)}}{
        \|\psi\|_{\dot{H}^\gamma(\R^n)}} = \varepsilon^{\frac {m+2}
        2\bigl( \frac n 2 - \gamma\bigr) -\frac 2 {\kappa -1} },
\]
and $\frac {m+2} 2\bigl( \frac n 2 - \gamma \bigr) -\frac 2 {\kappa -1}
> 0$ for $\gamma < \frac n 2 -\frac 4{(m+2)(\kappa-1)}$. Hence,
$\gamma < \frac n 2-\frac 4{(m+2)(\kappa-1)}$ implies that both the
norm of the data $(\varphi_\varepsilon, \psi_\varepsilon)$ and the
lifespan $T_{\ve}= \ve T$ of the solution $u_{\ve}$ go to zero as
$\ve\to 0$, where $T$ is the lifespan of the solution $u$.
\end{remark}

In case $ \kappa_\ast \le \kappa < \kappa_2$, as a supplement to
Theorem~\ref{thm1.1}, we consider the local existence and uniqueness
of solutions $u$ of problem \eqref{1.1} in the space $C([0, T];
\dot{H}^\gamma(\R^n)) \cap L^s((0, T); L^q (\R^n))$ for certain $s\neq
q$.

\begin{theorem} \label{thm1.5}
Let $n\geq2$, $F$ be above, $\gamma = \gamma(\kappa,m,n)$ be as
in\/ \textup{Theorem 1.1}, and suppose that $\kappa_\ast \le \kappa <
\kappa_2$. Then the unique solution $u$ of problem \eqref{1.1} also
belongs to the space $L^s((0, T);L^q(\R^n))$, where
\begin{align*}
  \frac 1 q&=\frac 1 {(m+2)(n-1)} \left( \frac 8{\kappa-1} - \frac m
        {\mu_\ast} \right) -\frac {n-1}{2(n+1)} \\
  \intertext{and}
  \frac 1 s&= \frac{(m+2)(n-1)}4 \left( \frac 1 2 -\frac 1 q\right)
  +\frac m {4 \mu_\ast}.
\end{align*}
Moreover, estimate \eqref{1.4} is satisfied.
\end{theorem}

\smallskip

If $n \ge 3$ or $n=2$, $m \ge 3$, then we find a number
$\gamma(\kappa, m, n)$ also for certain $\kappa$ in the range
$1<\kappa<\kappa_1$.

\begin{theorem} \label{thm1.6}
Let $n \ge 3$ or $n=2$, $m \ge 3$. Let $F$ be as above and
$\kappa_0 \le \kappa < \kappa_1$. In addition, let the exponent
$\gamma = \gamma(\kappa, m, n)$ in \eqref{1.1} be given by
\begin{equation}  \label{1.5}
\gamma(\kappa, m, n)=\frac {n+1} 4 -\frac{n+1}{4 \mu_\ast (m+2)}\cdot \frac{\mu_\ast
  (m+2)(n-1)+12 \mu_\ast +2m}{2n \kappa
  -(n+1)}-\frac m {2 \mu_\ast (m+2)}.
\end{equation}
Then problem \eqref{1.1} possesses a unique solution $u \in C([0, T];
\dot{H}^\gamma(\R^n)) \cap L^s( (0, T);L^q(\R^n))$ for some $T>0$,
where
\begin{align*}
  \frac 1 q &=\frac 1 {2n \kappa-(n+1)} \left(\frac{n-1} 2 +\frac
  6{m+2} +\frac m {\mu_\ast(m+2)} \right) \\
\intertext{and}
\frac 1 s &= \frac{(m+2)(n-1)}4 \left( \frac 1 2 -\frac 1 q\right)
+\frac m {4 \mu_\ast}.
\end{align*}
Moreover, estimate \eqref{1.4} is satisfied.
\end{theorem}


\begin{remark}
Other than for the wave equation when $m=0$ (see also
Remark~\ref{rem1.9} below), here $\gamma$ can be negative in certain
situations. In fact, $\gamma(\kappa,m,n)<0$ holds in the following
cases:

\begin{enumerate}[(i)]  
\item $\kappa_1<\kappa<\frac{35}{17}$ \ ($<\kappa_*$) if $n=2$, $m=1$
  and $\kappa_1<\kappa<\frac{13}7$ \ ($<\kappa_*$) if $n=2$, $m=2$
  (see Theorem~\ref{thm1.1}),
\item
  $\kappa_0<\kappa<\frac{\mu_*(\mu_*+2)(n+1)}{\mu_*(\mu_*-1)(n+1)-mn}$
  \ ($\leq \kappa_1$) if $n\geq3$ or $n=2$, $m\geq3$  (see Theorem~\ref{thm1.6}).
\end{enumerate}
\end{remark}

\begin{remark}  \label{rem1.3}
For initial data $(\varphi, \psi )$ belonging to $H^\gamma(\R^n)
\times H^{\gamma -\frac 2 {m+2}}(\R^n)$, where $\gamma \ge
\gamma(\kappa, m, n)$, Theorems~\ref{thm1.1}, \ref{thm1.5}, and
\ref{thm1.6} remain valid.
\end{remark}

\begin{remark}  \label{rem1.9}
For $m = 0$, \eqref{1.1} becomes
\begin{equation*}
\left\{ \enspace
\begin{aligned}
& \p_t^2 u - \Delta u =F(u) &&\text{in }   (0,T) \times \R^n, \\
& u (0, \cdot)=\varphi  \in \dot{H}^\gamma(\R^n), \\
 & \p_t u (0, \cdot)=\psi \in \dot{H}^{\gamma-1}(\R^n),
\end{aligned}
\right.
\end{equation*}
while the exponents $\kappa_\ast, \kappa_0, \kappa_1, \kappa_2$, and
$\kappa_3$ are
\[
\begin{aligned}
\kappa_\ast&=\frac {n+3}{n-1}, \quad
\kappa_2=\frac {(n+1)^2-6}{(n-1)^2-2}, \\
\kappa_1&=\frac{(n+1)^2}{(n-1)^2 +4} \quad \text{if } n \ge 3, \\
\kappa_0&=1+\frac3n, \quad
\kappa_3=\frac {n+1}{n-3} \quad \text{if } n \ge 4.
\end{aligned}
\]
For $n \ge 3$, $\gamma$ defined in \eqref{1.3} equals
\begin{equation}\label{1.6}
\gamma(\kappa,0,n)  =
\begin{cases}
\ds\frac {n+1} 4 -\frac{1}{\kappa-1} & \text{if } \kappa_1 < \kappa \le \kappa_*, \\
\ds\frac n 2 -\frac 2 {\kappa-1} & \text{if } \kappa\ge \kappa_*,
\end{cases}
\end{equation}
whereas, for $n \ge 4$, $\gamma$ defined in \eqref{1.5} equals
\begin{equation}\label{1.7}
\gamma(\kappa,0,n)=\frac {n+1} 4 -\frac{(n+1)(n+5)}4\,\frac1{2n \kappa -(n+1)} .
 \end{equation}
Note that the numbers in \eqref{1.6} and \eqref{1.7} are exactly those
in (2.1) and (2.5) of \cite{Lind4}. In that paper, \cite{Lind4}, the
local existence problem for minimal regularity solutions of the
semilinear wave equation was systematically studied. The results were
achieved by establishing Strichartz-type estimates for the linear wave
operator $\p_t^2 - \Delta$. Under certain restrictions on the
nonlinearity $F(u, \na u)$, for the more general semilinear wave
equation
\begin{equation*}
\left\{ \enspace
\begin{aligned}
&\partial_t^2 u-\Delta u = F(u, \na u),  \\
&u(0,x)=\varphi(x), \quad \partial_{t} u(0,x)=\psi(x),\\
\end{aligned}
\right.
\end{equation*}
many remarkable results on the ill-posedness or well-posedness problem
on the local existence of low regularity solutions have been obtained,
see \cite{Kap1, Lind1, Lind4, Ponce, Tat1, Stru} and the reference
therein.

\end{remark}

\begin{remark}\label{rem1.10}
There are some essential differences between degenerate hyperbolic
equations and \linebreak strictly hyperbolic equations. Amongst
others, the symmetry group is smaller (see \cite{Lup1}) and there is a
loss of regularity for the linear Cauchy problem (see
e.g.~\cite{Dreher,Taniguchi}). Therefore, as compared to the
semilinear wave equation, a more delicate analysis is required when
one studies minimal regularity results for the semilinear generalized
Tricomi equation in the degenerate hyperbolic region.
\end{remark}

\smallskip

The Tricomi equation (i.e., Eq.~\eqref{1.1} for $n=1$, $m=1$) were
first studied by Tricomi \cite{Tricomi1} who initiated work on
boundary value problems for linear partial differential operators of
mixed elliptic-hyperbolic type. So far, these equations have been
extensively studied in bounded domains under suitable boundary
conditions and several applications to transonic flow problems were
given (see \cite{Bers, Germain, Tricomi1, Morawetz} and the references
therein). Conservation laws for equations of mixed type were derived
by Lupo and Payne \cite{Lup1, Lup2}. In \cite{RWY3}, we established the
local solvability for low regularity solutions of the semilinear
equation $ \p_t^2 u - t^m \Delta u = F(u)$, where $n\ge 2$, $m \in \N$
is odd, in the domain $(-T, T) \times \R^n$ for some $T>0$.
In \cite{Gelf, Yag1, Yag3}, fundamental solutions for the linear
Tricomi operator and the linear generalized Tricomi operator have been
explicitly computed. In case $n=2$ and $m=1$, Beals \cite{Beals}
obtained the local existence of the solution $u$ of the equation
$\p_t^2u-t\Delta u=F(u)$ with initial data of $H^s$-regularity, where
$s>n/2$. For the equation $\p_t^2 u- t^m \Delta u = a(t) F(u)$, where
$n \ge 2$, $m \in\mathbb N$ is even, and both $a$ and $F$ are of power
type, Yadgjian \cite{Yag2} obtained global existence and uniqueness
for small data solutions provided the solution $v$ of the linear
problem $\p_t^2v - t^m \Delta v = 0$ fulfills $t^\beta v\in C([0,
  \infty); L^q(\R^n))$ for certain $\beta,\,q$ depending on $n$, $m$,
  and the powers occurring in $a$ and $F$.
In \cite{RWY1, RWY2}, for the semilinear generalized Tricomi equation
$ \p_t^2 u -t^m \Delta u =F(u)$ with initial data of a special
structure, i.e., homogeneous of degree $0$ or piecewise smooth along a
hyperplane, we obtained local existence and uniqueness via
establishing $L^\infty $ estimates on the solutions $v$ of the linear
equation $\p_t^2 v -t^m \Delta v=g$. Note that when the nonlinear term
$F(u)$ is of power type, for higher and higher powers of $\kappa$,
these $L^\infty$ estimates are basically required to guarantee
existence. In this paper, where the initial data in
$\dot{H}^\gamma(\R^n)$ is of no special structure and $\gamma$ is
minimal to guarantee local well-posedness of problem (1.1),
the arguments of
\cite{RWY1, RWY2} fail. Inspired by the methods in \cite{Lind4},
however, we are able to overcome the technical difficulties related to
degeneracy and low regularity and eventually obtain the local
well-posedness of problem \eqref{1.1}.

\medskip

We first study the linear problem
\begin{equation}\label{linear_problem}
\left\{ \enspace
\begin{aligned}
& \p_t^2 u -t^m \Delta u = f(t,x) &&\text{in $(0,T)\times \R^n$,} \\
& u(0,\cdot)=\varphi(x),  \quad  \p_t u(0,\cdot) =\psi(x)
\end{aligned}
\right.
\end{equation}
and establish Strichartz-type estimates of the form
\begin{equation}\label{strichartz}
  \|u\|_{C_t^0\dot{H}_x^\gamma(S_T)} +
  \|u\|_{L_t^sL_x^q(S_T)} \leq C \left( \|\varphi\|_{\dot{H}^\gamma} +
  \|\psi\|_{\dot{H}^{\gamma-\frac2{m+2}}} +
  \|f\|_{L_t^rL_x^p(S_T)}\right)
\end{equation}
for certain $s,\,q,\,r,\,p$ (for details see below) and some constant
$C=C(T,\gamma,s,q,r,p)>0$, where $S_T=(0,T)\times\R^n$. Note that, by
scaling, a necessary condition for this estimate in case $T=\infty$ to
hold is
\begin{equation}  \label{SI}
\frac{\left(m+2\right)n}{2} \left( \frac 1 p -\frac 1 q \right) + \frac 1 r- \frac
1 s =2.
\end{equation}

\medskip

In doing so, in Section~\ref{sec2}, we introduce certain Fourier
integral operators $W$ ($=W^0$) and $W^\alpha$ for $\alpha\in\mathbb
C$. These operators depend on a parameter $\mu\ge 2$, introduced in
\eqref{2.17}, which plays an auxiliary role for the linear problems
and agrees with the homogeneous dimension $\mu_*$ when applied to the
semilinear problems. Along with the operators $W$ and $W^\alpha$ we
also consider their parts $W_j$ and $W_j^\alpha$, respectively,
resulting from a dyadic decomposition of frequency space. Continuity
of the operators $W_j$ and $W_j^\alpha$ between function spaces which
holds uniformly in $j$ ultimately provides linear estimates on the
solutions $u$ of Eq.~\eqref{linear_problem}.


In Section~\ref{sec3}, we prove boundedness of the operators
$W^\alpha_j$ from $L^r_t L^p_x(\R^{1+n}_+)$ to $ L^{r'}_t
L^{p'}_x(\R^{1+n}_+)$ (see Theorem~\ref{thm3.1}) and from $L^r_t
L^p_x(\R^{1+n}_+)$ to $L^\infty_t L^2_x(\R^{1+n}_+)$ (see
Theorem~\ref{thm3.4}), where $\mu$ has to satisfy the lower bound $
\mu \ge \max\{2,m/2\}$. Combining Theorem~\ref{thm3.1} and Stein's
analytic interpolation theorem, we show boundedness of the operators
$W^\alpha_j$ from $L^q(\R^{1+n}_+)$ to $ L^{p_0}(\R^{1+n}_+)$, where $
q_0 \le q \le \infty$ (see Theorem~\ref{thm3.6}). Through an
additional dyadic decomposition now with respect to the time variable
$t$, using Theorems~\ref{thm3.1} and \ref{thm3.6} together with
interpolation, we prove boundedness of the operators $W_j$ from $L^r_t
L^p_x((0, T) \times \R^n)$ to $ L^s_t L^q_x((0, T) \times \R^n)$ for
any $T>0$ (see Theorems~\ref{thm3.7} and \ref{thm3.8}), where $\mu$
has to satisfy the new lower bounds $\mu \ge \mu_\ast$
(Theorem~\ref{thm3.7}) and $\mu \ge \max\{2, mn/2 \}$
(Theorem~\ref{thm3.8}), respectively.

\medskip

In the sequel, we shall use the following notation:
\[
\frac 1 {p_0}= \frac 1 2 +\frac{2\mu-m}{\mu(2\mu_\ast -m)}, \quad
\frac 1 { p_1}=\frac 1 2+ \frac { 2\mu -m}{\mu (m+2)(n-1)}, \quad
\frac 1 {p_2}=\frac 2 {p_0} -\frac 1 {p_1}.
\]
Note that
\[
1< p_1\leq p_0\leq p_2\leq 2 \qquad \text{if $n\geq3$ or $n=2$, $m\geq2$,}
\]
while $1\leq p_1$ in case of $n=2$ and $m=1$ requires $\mu=2$ (and then
$p_1=1$). For $1\leq p\leq 2$, $p'$ denotes the conjugate exponent of
$p$ defined by $\ds\frac1p+\frac1{p'}=1$. Further, $q_\ell$ denotes
$p_\ell'$ for $\ell=0,1,2$, while $q_0^*$ equals $q_0$ when
$\mu=\mu_*$ (see Remark~\ref{rem4.2}). We often abbreviate function
spaces $C_t^0\dot{H}^\gamma_x(S_T)=C([0,T]; \dot{H}^\gamma(\R^n))$,
$L^r_tL^p_x(S_T)=L^r((0,T); L^p(\R^n))$, and $A\lesssim B$ means that
$A\le CB$ holds for some generic constant $C>0$.

\medskip

The paper is organized as follows: In Section~\ref{sec2}, we define a
class of Fourier integral operators associated with the linear
generalized Tricomi operator $\p_t^2-t^m\Delta$ in $\R_+ \times \R^n$.
Then, in Section~\ref{sec3}, we establish a series of mixed-norm
space-time estimates for those Fourier integral operators.  These
estimates are applied, in Section~\ref{sec4}, to obtain
Strichartz-type estimates for the solutions of the linear generalized
Tricomi equation which in turn, in Section~\ref{sec5}, allow us to
prove the local existence and uniqueness results for problem
\eqref{1.1}.


\section{Some preliminaries}\label{sec2}

In this section, we first recall an explicit formula for the solution
of the linear generalized Tricomi equation obtained in
\cite{Taniguchi} and then apply it to define a class of Fourier
integral operators which will play a key role in proving our main
results.

Consider the linear generalized Tricomi equation
\begin{equation}\label{2.1}
\left\{ \enspace
\begin{aligned}
& \p_t^2 u -t^m \Delta u =f(t,x)   && \text{in }    \R_{+} \times \R^n, \\
& u (0, \cdot)=\varphi,  \quad \p_t u(0, \cdot)=\psi.
\end{aligned}
\right.
\end{equation}
Its solution $u$ can be written as $u=v + w$, where $v$ solves the
homogeneous equation
\begin{equation} \label{2.2}
\left\{ \enspace
\begin{aligned}
& \p_t^2 v -t^m \Delta v =0    && \text{in }    \R_{+} \times \R^n, \\
& v (0, \cdot)=\varphi, \quad \p_t v (0, \cdot)=\psi
\end{aligned}
\right.
\end{equation}
and $w$ solves the inhomogeneous equation with zero initial data
\begin{equation}\label{2.3}
\left\{ \enspace
\begin{aligned}
& \p_t^2 w -t^m \Delta w =f(t,x)  && \text{in }      \R_{+} \times \R^n, \\
& w (0, \cdot)= \p_t w (0, \cdot)=0.
\end{aligned}
\right.
\end{equation}

Recall that (see \cite{Taniguchi} or \cite{Yag2}) the solutions $v$ and $w$ of
problems \eqref{2.2} and \eqref{2.3} can be expressed as
\[ v(t,x) = V_0(t, D_x) \varphi(x) +V_1(t, D_x) \psi(x)
\] 
and
\begin{align} \label{2.5}
w(t,x)= \int_0^t \left( V_1(t, D_x) V_0(\tau, D_x) -   V_0(t, D_x) V_1(\tau,
D_x) \right) f(\tau, x)\, d\tau,
\end{align}
where the symbols $V_j(t,\xi)$ ($j=0, 1$) of the Fourier
integral operators $V_j(t, D_x)$ are
\begin{equation}\label{2.6}
\left\{ \enspace
\begin{aligned}
V_0(t,\xi)&=e^{-z/2}\,\Phi\left(\f{m}{2(m+2)},\f{m}{m+2};z\right),\\
V_1(t,\xi)&=te^{-z/2}\,\Phi\left(\f{m+4}{2(m+2)},\f{m+4}{m+2};z\right)
\end{aligned}
\right.
\end{equation}
with $z=2i\phi(t)|\xi|$ and $\phi(t)= \left(2/(m+2)\right)
t^{(m+2)/2}$. Here, $\Phi(a,c;z)$ is the confluent hypergeometric
function which is an analytic function of $z$.  Recall (see
\cite[page254]{Erd1}) that
\begin{equation} \label{2.7}
\frac{d^n}{dz^n} \Phi(a,c;z)=\frac{(a)_n}{(c)_n} \Phi(a+n, c+n;z),
\end{equation}
where $(a)_0 = 1$, $(a)_n= a(a+1) \ldots (a+n-1)$. In addition (see
\cite[(3.5)-(3.7)]{Yag2}), for $0< \operatorname{arg(z)}<\pi$, one has that
\begin{equation} \label{2.8}
 e^{-z/2}\,\Phi(a,c;z)
= \frac{\Gamma(c)}{\Gamma(a)}\, e^{z/2}H_+(a,c;z)
+\frac{\Gamma(c)}{\Gamma(c-a)}\, e^{-z/2}H_-(a,c;z),
\end{equation}
where
\begin{equation*}
\begin{aligned}
H_+(a,c;
z)&=\f{e^{-i\pi(c-a)}}{e^{i\pi(c-a)}-e^{-i\pi(c-a)}}\f{1}{\G(c-a)}\,z^{a-c}
\int_{\infty}^{(0+)}e^{-\th}\th^{c-a-1}\left(1-\f{\th}{z}\right)^{a-1}d\th,\\
H_-(a,c; z)&=\f{1}{e^{i\pi a}-e^{-i\pi a}}\f{1}{\G(a)}\,z^{-a}
\int_{\infty}^{(0+)}e^{-\th}\th^{a-1}\left(1+\f{\th}{z}\right)^{c-a-1}d\th.\\
\end{aligned}
\end{equation*}
Moreover, it holds that
\begin{equation}\label{2.9}
\begin{aligned}
  \left| \partial_\xi^\beta\left( H_{+}(a,c;2i\phi(t)|\xi|) \right)
  \right| \lesssim & (\phi(t)|\xi|)^{a-c} (1+ |\xi|)^{ -|\beta|} &&
  \text{if } \phi(t)|\xi|\geq 1, \\
  \left|\partial_\xi^\beta \left(
  H_{-}(a,c;2i\phi(t)|\xi|)\right)\right| \lesssim &
  (\phi(t)|\xi|)^{-a}(1+|\xi|)^{-|\beta|} && \text{if }
  \phi(t)|\xi|\geq 1.
\end{aligned}
\end{equation}

Choose $\eta \in C^\infty_c(\R_+)$ such that $0 \le \eta \le 1$ with
$\eta(r) = 1$ if $r\le 1$ and $\eta(r)= 0$ if $r\geq2$. Then from
\eqref{2.6} and \eqref{2.8}, we can write
\begin{equation}\label{2.11}
 V_0(t, D_x) \varphi(x)= \int_{\R^n} e^{i (x \cdot \xi -
   \phi(t) |\xi|)}b_1(t, \xi) \hat{\varphi}(\xi)\, \dbar\xi  +
 \int_{\R^n} e^{i (x \cdot \xi + \phi(t) |\xi|)} b_2(t,
 \xi) \hat{\varphi}(\xi) \ \dbar\xi
\end{equation}
and
\begin{equation}\label{2.12}
  V_1(t, D_x) \psi(x)= \int_{\R^n} e^{i (x \cdot \xi - \phi(t) |\xi|)}
  b_3(t, \xi) \hat{\psi}(\xi) \, \dbar\xi + \int_{\R^n} e^{i (x
    \cdot \xi + \phi(t) |\xi|)} b_4(t, \xi) \hat{\psi}(\xi) \,
  \dbar\xi,
\end{equation}
where
\begin{align*}
b_1(t,\xi)&=\eta(\phi(t)
|\xi|)\Phi\!\left(\f{m}{2(m+2)},\f{m}{m+2};z\right)  +\bigl( 1-\eta(\phi(t)
|\xi|)\bigr) H_-\left(\f{m}{2(m+2)},\f{m}{m+2};z\right), \\
b_2(t,\xi)&=\bigl( 1-\eta(\phi(t) |\xi|)\bigr)
H_+\left(\f{m}{2(m+2)},\f{m}{m+2};z\right), \\
\intertext{and}
b_3(t,\xi)&= t\eta(\phi(t)
|\xi|)\Phi\!\left(\f{m+4}{2(m+2)},\f{m+4}{m+2};z\right) + t\bigl(
1-\eta(\phi(t) |\xi|)\bigr)
H_-\left(\f{m+4}{2(m+2)},\f{m+4}{m+2};z\right), \\
b_4(t,\xi)&= t
\bigl( 1-\eta(\phi(t) |\xi|)\bigr)
H_+\left(\f{m+4}{2(m+2)},\f{m+4}{m+2};z\right).
\end{align*}
We can also write
\begin{multline}\label{2.13}
\int_0^t V_0(t, D_x) V_1(\tau, D_x) f(\tau,x) \, d\tau = \int_0^t
\int_{\R^n} e^{i \left(x \cdot \xi + (\phi(t)+\phi(\tau))|\xi|\right)}
b_2(t,\xi) b_4(\tau, \xi ) \hat{f}(\tau, \xi) \ \dbar\xi d\tau
\\
\begin{aligned}
& + \int_0^t \int_{\R^n} e^{i \left(x \cdot \xi +
(\phi(t)-\phi(\tau))|\xi|\right)} b_2(t,\xi) b_3(\tau, \xi )  \hat{f}(\tau,
\xi)  \ \dbar\xi d\tau\\
& + \int_0^t \int_{\R^n} e^{i \left(x \cdot \xi -
(\phi(t)+\phi(\tau))|\xi|\right)}b_1(t,\xi) b_3(\tau, \xi )   \hat{f}(\tau,
\xi)  \ \dbar\xi d\tau \\
& + \int_0^t \int_{\R^n} e^{i \left(x \cdot \xi -
(\phi(t)-\phi(\tau))|\xi|\right)}b_1(t,\xi) b_4(\tau, \xi )   \hat{f}(\tau,
\xi)  \ \dbar\xi d\tau
\end{aligned}
\end{multline}
and
\begin{multline} \label{2.14}
\int_0^t  V_1(t, D_x) V_0(\tau, D_x)  f(\tau,x) d\tau
= \int_0^t \int_{\R^n} e^{i \left(x \cdot \xi +
(\phi(t)+\phi(\tau))|\xi|\right)}   b_4(t,\xi) b_2(\tau, \xi ) \hat{f}(\tau,
\xi) \ \dbar\xi d\tau\\
\begin{aligned}
&  + \int_0^t \int_{\R^n} e^{i \left(x \cdot \xi -
(\phi(t)-\phi(\tau))|\xi|\right)} b_3(t,\xi) b_2(\tau, \xi )     \hat{f}(\tau,
\xi) \ \dbar\xi d\tau\\
&  + \int_0^t \int_{\R^n} e^{i \left(x \cdot \xi -
(\phi(t)+\phi(\tau))|\xi|\right)} b_3(t,\xi)b_1(\tau, \xi )    \hat{f}(\tau,
\xi) \ \dbar\xi d\tau \\
& + \int_0^t \int_{\R^n} e^{i \left(x \cdot \xi +
(\phi(t)-\phi(\tau))|\xi|\right)} b_4(t,\xi)b_1(\tau, \xi )   \hat{f}(\tau,
\xi) \ \dbar\xi d\tau,
\end{aligned}
\end{multline}
where $\hat{f}(\tau, \xi)$ is the Fourier transform of $f(\tau,x)$
with respect to the variable $x$ and $\dbar\xi=(2\pi)^{-n}\,d\xi$.

In view of the analyticity of $\Phi(a,c;z)$ with respect to the
variable $z$, identity \eqref{2.7}, and estimates~\eqref{2.9}, we have
that, for $(t,\xi) \in \R_+^{1+n}$,
\begin{align} \label{2.15}
 \bigl|\partial_\xi^\beta b_\ell(t, \xi)\bigr| \lesssim  \bigl( 1+ \phi(t)
|\xi|\bigr)^{-\frac{m}{2(m+2)}} |\xi|^{-|\beta|}, \quad \ell=1,2,
\end{align}
and
\begin{align}\label{2.16}
 \bigl |\partial_\xi^\beta b_\ell(t, \xi)\bigr| \lesssim t \bigl( 1+
 \phi(t) |\xi|\bigr)^{-\frac{m+4}{2(m+2)}} |\xi|^{-|\beta|}, \quad
 \ell=3,4.
\end{align}
Thus, for $\ell=1,2$, $k=3,4$, $\mu \ge 2$, $t,\,  \tau >0$,
and $\xi \in \R^n$, one has from \eqref{2.15} and \eqref{2.16} that
\begin{equation}\label{2.17}
\begin{aligned}
\bigl|\partial_\xi^\beta \bigl(b_k(t, \xi) b_\ell(\tau,
\xi) \bigr) \bigr| & \lesssim t \bigl( 1+ \phi(t) |\xi|
\bigr)^{-\frac{m+4}{2(m+2)}} \bigl( 1+ \phi(\tau) |\xi|
\bigr)^{-\frac{m}{2(m+2)}} |\xi|^{-|\beta|} \\ & \lesssim
\bigl( 1+ \phi(t) |\xi| \bigr)^{-\frac{m}{2(m+2)}} \bigl( 1+ \phi(\tau)
|\xi| \bigr)^{-\frac{m}{2(m+2)}}
|\xi|^{-\frac{2}{m+2}-|\beta|}\\ & \lesssim \bigl( 1+
\left|\phi(t)-\phi(\tau)\right| |\xi| \bigr)^{-\frac{m}{\mu(m+2)}}
|\xi|^{-\frac{2}{m+2}-|\beta|}.
\end{aligned}
\end{equation}
Furthermore, estimates \eqref{2.15}-\eqref{2.17} yield that, for
$\ell=1,2$, $k=3,4$ or $\ell=3,4$, $k=1,2$ and for $\mu \ge 2$,
$t,\,s>0$, and $\xi \in \R^n$, one has
\begin{multline}\label{2.18}
\left| \partial_\xi^\beta \Bigl( \int_t^\infty \overline{b_{\ell}(\tau, \xi) b_{k}(t, \xi) } \ \p_\tau\bigl( b_{\ell}(\tau, \xi) b_{k}(s, \xi) \bigr) d\tau \Bigr) \right|
\\ \lesssim \bigl( 1+ \big|\phi(t)-\phi(s)\big| |\xi|
\bigr)^{-\frac{m}{\mu(m+2)}} |\xi|^{-\frac{4}{m+2}-|\beta|}
\end{multline}
and
\begin{multline}\label{2.19}
\left| \partial_\xi^\beta \Bigl( \int_s^\infty \overline{b_{\ell}(\tau, \xi) b_{k}(t, \xi) } \ \p_\tau\bigl( b_{\ell}(\tau, \xi) b_{k}(s, \xi) \bigr) d\tau \Bigr) \right|
\\ \lesssim \Bigl( 1+ \big|\phi(t)-\phi(s)\big| |\xi|
\Bigr)^{-\frac{m}{\mu(m+2)}} |\xi|^{-\frac{4}{m+2}-|\beta|}.
\end{multline}

In order to study the function $w$ in \eqref{2.5}, in view of
\eqref{2.13}, \eqref{2.14} and \eqref{2.17}-\eqref{2.19}, it suffices
to consider, for a given $\mu \ge 2$, the Fourier integral operator
$W$,
\begin{align} \label{2.20}
  Wf(t,x)= \int_0^t \int_{\R^n} e^{i \left(x \cdot \xi +
    (\phi(t)-\phi(s)) |\xi|\right)} b(t, s, \xi) \hat{f}(s, \xi)
  \ \dbar\xi ds,
\end{align}
where  $b \in C^\infty(\R_+ \times \R_+ \times \R^n)$ satisfies

\noindent
(i) \ for $t,\, s >0$ and $\xi \in \R^n$,
\begin{equation} \label{2.21}
\bigl|\partial_\xi^\beta b(t, s, \xi)\bigr| \lesssim \left( 1+
\left|\phi(t)-\phi(s)\right| |\xi| \right)^{-\frac{m}{\mu(m+2)}}
|\xi|^{-\frac{2}{m+2}-|\beta|},
\end{equation}

\noindent
(ii) \ for $t,\, s > 0$  and $\xi \in \R^n$,
\begin{equation}\label{2.22}
\Bigl| \partial_\xi^\beta \Bigl( \int_t^\infty \overline{b(\tau, t,
  \xi) } \ \p_\tau b(\tau, s, \xi) d\tau \Bigr) \Bigr| \lesssim \left(
1+ \left|\phi(t)-\phi(s)\right| |\xi| \right)^{-\frac{m}{\mu(m+2)}}
|\xi|^{-\frac{4}{m+2}-|\beta|}
\end{equation}
and
\begin{equation}\label{2.23}
\Bigl| \partial_\xi^\beta \Bigl( \int_s^\infty \overline{b(\tau, t,
  \xi) } \ \p_\tau b(\tau, s, \xi) \, d\tau \Bigr) \Bigr| \lesssim
\left( 1+ \left|\phi(t)-\phi(s)\right| |\xi| \right)^{-\frac{m}{\mu(m+2)}}
|\xi|^{-\frac{4}{m+2}-|\beta|}.
\end{equation}

Let $\Theta \in C^\infty_c(\R_+)$ satisfy
$\operatorname{supp}\Theta\subseteq[1/2,2]$ and
\[
\sum_{j=-\infty}^{\infty}\Theta(t/2^j)=1  \quad \text{for }  t>0.
\]
Then, as in \cite{Lind4}, for $j\in \mathbb Z$ and $\alpha\in
\mathbb{C}$, we define dyadic operators $W_j$ and $W_j^\alpha$,
\[
  W_j f (t,x)= \int_0^t \int_{\R^n} e^{i \left(x \cdot \xi +
    (\phi(t)-\phi(s)) |\xi|\right)} b_j(t, s, \xi) \hat{f}(s, \xi)
  \ \dbar\xi ds
\]
and
\begin{equation}\label{2.25}
  W_j^\alpha f(t,x)= \int_0^t \int_{\R^n} e^{i \left(x \cdot \xi +
    (\phi(t)-\phi(s)) |\xi|\right)} b_j(t, s, \xi) \hat{f}(s, \xi)\,
  \frac {\dbar\xi}{|\xi|^\alpha} ds,
\end{equation}
where $b_j(t, s, \xi) = \Theta(|\xi|/2^j) b(t, s, \xi) $; here $b
\in C^\infty(\R_+ \times \R_+ \times \R^n)$ satisfies the estimates
\eqref{2.21}-\eqref{2.23}.

\smallskip

Littlewood-Paley theory gives us a relationship between $Wf$ and
$W_jf$ ($= W^0_j f$), which will play an important role in our
arguments in Section~\ref{sec4}.

\begin{proposition} \label{prop2.1}
Let $n \ge 2$. For $1 < p \le 2$, $1\le r \le 2$, $2 \le q <\infty$,
and $2\le s \le \infty$, let
\begin{equation} \label{2.26}
\|W_j f\|_{L^s_t L^q_x}  \lesssim \|f\|_{L^r_t L^p_x}
\end{equation}
hold uniformly in $j$. Then
\[\|W f\|_{L^s_t L^q_x} \lesssim  \|f\|_{L^r_t L^p_x}.\]
\end{proposition}

\begin{proof}
This is actually an application of Lemma~3.8 of \cite{Lind4}.  For the
sake of completeness, we give the proof here. By Littlewood-Paley
theory (see, e.\,g., \cite{Stein}), for any $1<\rho<\infty$,
\begin{equation*}
 \|Wf(t, \cdot) \|_{L^\rho(\R^n)} \lesssim \Bigl\| \Bigl(\sum
 \limits_{j=-\infty}^\infty |W_jf(t, \cdot)|^2
 \Bigr)^{1/2}\Bigr\|_{L^\rho(\R^n)} \lesssim \|Wf(t, \cdot)
 \|_{L^\rho(\R^n)}.
\end{equation*}
Together with the Minkowski inequality, this yields
\begin{equation}\label{2.27}
  \|Wf \|_{L^s_t L^q_x} \lesssim \Bigl(\sum \limits_{j=-\infty}^\infty
  \|W_jf\|_{L^s_t L^q_x}^2 \Bigr)^{1/2}
\end{equation}
and
\begin{equation} \label{2.28}
  \Bigl(\sum \limits_{j=-\infty}^\infty \|W_jf\|_{L^r_t L^p_x}^2
  \Bigr)^{1/2} \lesssim \|Wf \|_{L^r_t L^p_x}.
\end{equation}
Notice that
\[
  f =\sum \limits_{k=-\infty}^\infty f_k,
\]
where $f_k(\tau, x)= \Theta\left(\tau/2^k\right) f(\tau,
x)$. Therefore, for some $M_0\in\N$,
\begin{align*}
\|W f\|_{L^s_t L^q_x}^2 &\lesssim \sum \limits_{j=-\infty}^\infty
\|W_jf\|_{L^s_t L^q_x}^2 && \text{(by \eqref{2.27})} \\
&= \sum \limits_{j=-\infty}^\infty \Bigl\|W_j\Bigl( \ \sum
\limits_{|j-k|\le M_0}f_k\Bigr)\Bigr\|_{L^s_t L^q_x}^2 && \text{(due
  to the compact support of } \Theta) \\
&\lesssim \sum \limits_{j=-\infty}^\infty \Bigl( \,
\sum\limits_{|j-k|\le M_0} \|W_j f_k\|_{L^s_t L^q_x} \Bigr)^2 &&
\text{(by Minkowski inequality)} \\
&\lesssim \sum \limits_{j=-\infty}^\infty \sum \limits_{|j-k|\le M_0 }
\|f_k\|_{L^r_t L^p_x}^2 && \text{(by \eqref{2.26})}\\
&\lesssim \sum \limits_{j=-\infty}^\infty \|f_j\|_{L^r_t L^p_x}^2
\lesssim \|f\|_{L^r_t L^p_x}^2. && \text{(by \eqref{2.28})},
\end{align*}
which completes the proof of Proposition~\ref{prop2.1}.
\end{proof}



\section{Mixed-norm estimates for a class of Fourier integral operators}\label{sec3}

In this section, for $j \in \mathbb Z$, $\alpha \in \mathbb C$, and
$\mu \ge 2$, we shall study mixed norm estimates for the class of
Fourier integral operators $W^\alpha_j$ defined in \eqref{2.25}.

We start by considering the boundedness of the operator $W^\alpha_j$
from $L^r_t L^p_x$ to $ L^{r'}_t L^{p'}_x$, where $1< r,\, p \le
2$. We denote $\lambda_j=2^j$. \textit{All the following estimates
  hold uniformly in $j$.}

\begin{theorem}\label{thm3.1}
Let $n \ge 2$ and $\mu \ge \max\{2,m/2\}$. Then\/\textup{:}

\noindent
\textup{(i)} \ For $\max\{p_1,1\} <p \le 2$ and
\begin{equation} \label{3.1}
 \frac 1 r = 1-\frac m {4\mu} -\frac {(m+2)(n-1)} 4 \left( \frac 1 p
 -\frac 1 2 \right),
\end{equation}
we have that
\begin{equation} \label{3.2}
  \|W_j^\alpha f\|_{L^{r'}_t L^{p'}_x(\R^{1+n}_+)} \lesssim
  \lambda_j^{\left(\frac{1}{p}-\frac{1}{2}\right)(n+1)-\frac m {\mu (m+2)}-\frac
    2 {m+2}-\textup{Re }\alpha} \|f\|_{L^r_t L^p_x(\R^{1+n}_+)}.
\end{equation}
Consequently,
\begin{equation} \label{3.3}
  \|W_j^\alpha f\|_{L^{r'}_t L^{p'}_x(\R^{1+n}_+)} \lesssim
  \|f\|_{L^r_t L^p_x(\R^{1+n}_+)} \quad \textup{if } \ \operatorname{Re}
  \alpha=\left(\frac{1}{p}-\frac{1}{2}\right)(n+1)-\frac m {\mu
    (m+2)}-\frac 2 {m+2}.
\end{equation}

\noindent
\textup{(ii)} \  For  $ p_1>1$ and $ 1<p<p_1$, we have that
\begin{equation} \label{3.4}
  \|W_j^\alpha f\|_{L^{2}_t L^{p'}_x(\R^{1+n}_+)} \lesssim
  \lambda_j^{n\left(\frac{2}{p}-1\right)-\frac{4}{m+2}-\operatorname{Re}\alpha}
  \|f\|_{L^2_t L^p_x(\R^{1+n}_+)}.
\end{equation}
In particular,
\begin{equation} \label{3.5}
  \|W_j^\alpha f\|_{L^{2}_t L^{p'}_x(\R^{1+n}_+)} \lesssim
  \|f\|_{L^2_t L^p_x(\R^{1+n}_+)} \quad \textup{if } \
  \operatorname{Re}\alpha=n\left(\frac{2}{p}-1\right)-\frac{4}{m+2}.
\end{equation}
\end{theorem}

To prove Theorem~\ref{thm3.1}, for fixed $t,\,\tau >0$, we first
consider the operator $B_j^\alpha$,
\begin{equation*}
  B_j^\alpha f (t, \tau, x) = \int_{\R^n} e^{i \left( x \cdot \xi +
    (\phi(t)-\phi(\tau))|\xi|\right)} \ b_j(t, \tau, \xi)
  \hat{f}(\tau, \xi)\, \frac{\dbar\xi}{|\xi|^\alpha}.
\end{equation*}

\begin{lemma}\label{lem3.2}
Let $n \ge 2 $ and $1\le p \le 2$. Then, for $t,\tau> 0$,
\begin{multline}  \label{3.6}
  \|B_j^\alpha f (t, \tau, \cdot)\|_{L^{p'}(\R^n)} \lesssim
  \lambda_j^{\left(\frac{1}{p}-\frac{1}{2}\right)(n+1)-\frac m
    {\mu(m+2)}-\frac 2 {m+2}-\operatorname{Re}\alpha} \\ \times \bigl(
  \lambda_j^{-\frac{2}{m+2}} + |t-\tau|\bigr)^{-\left(m+2\right)\left(
    \frac{1}{p}-\frac{1}{2}\right)\frac{n-1}{2} - \frac{m}{2\mu}}
  \|f(\tau, \cdot)\|_{L^p(\R^n)}.
\end{multline}
\end{lemma}
\begin{proof}
Denote
\begin{equation} \label{3.7}
  K_j^\alpha(t,\tau, x, y)= \int_{\R^n} e^{i \left( (x-y) \cdot \xi +
    (\phi(t)-\phi(\tau))|\xi|\right)} \ b_j(t, \tau, \xi)\,
  \frac{\dbar\xi}{|\xi|^\alpha}.
\end{equation}
Then $B_j^\alpha f$ can be written as
\begin{equation*}
  B_j^\alpha f (t, \tau, x)=\int_{\R^n} K_j^\alpha(t,\tau, x, y)
  f(\tau, y)\, dy.
\end{equation*}
Since $\operatorname{supp}_\xi b_j \subseteq \{\xi \in \R^n\mid
\lambda_j/2\le |\xi| \le 2 \lambda_j \}$, we have from \eqref{2.21} that
\begin{align} \label{3.8}
|\partial_\xi^\beta b_j(t, \tau, \xi)| \lesssim
\lambda_j^{-\frac{m}{\mu(m+2)}-\frac 2 {m+2}-|\beta|} \bigl(
\lambda_j^{-\frac{2}{m+2}}+| t-\tau | \bigr)^{-\frac{m}{2\mu}}.
\end{align}
We now apply \eqref{3.8} to derive estimate \eqref{3.6} by
Plancherel's theorem when $p=2$ and by the stationary phase method
when $p=1$. By interpolation, we then obtain \eqref{3.6} for $1<p<2.$

Indeed, it follows from Plancherel's theorem that
\begin{align}
  \|B_j^\alpha f (t, \tau, \cdot)\|_{L^{2}_x(\R^n)} &=
  \bigl\|e^{i(\phi(t)-\phi(\tau))|\xi|} b_j(t, \tau, \xi)
  \hat{f}(\tau, \xi) |\xi|^{-\alpha}\bigr\|_{L^{2}_{\xi}(\R^n)} \notag
  \\ & \lesssim \lambda_j^{-\frac{m}{\mu(m+2)}-\frac 2
    {m+2}-\operatorname{Re} \alpha}\ \bigl(
  \lambda_j^{-\frac{2}{m+2}}+| t-\tau | \bigr)^{-\frac{m}{2\mu}}
  \ \|f(\tau, \cdot)\|_{L^2(\R^n)}. \label{3.9}
\end{align}
On the other hand, by the stationary phase method (see
e.g.~\cite[Lemma 7.2.4]{Sog}), we have that, for any $N \ge 0$,
\begin{align}
  |K_j^\alpha(t, \tau, x, y)| & \lesssim \lambda_j^n \bigl(
  1+\left|\phi(t)-\phi(\tau)\right|\lambda_j \bigr)^{-\frac{n-1}{2}}
  \bigl( \lambda_j^{-\frac{2}{m+2}}+| t-\tau |
  \bigr)^{-\frac{m}{2\mu}} \notag \\
  & \quad\qquad \times \lambda_j^{-\frac{m}{\mu(m+2)}-\frac 2
    {m+2}-\operatorname{Re} \alpha}\ \Bigl( 1+ \lambda_j\bigl||x-y|
  -|\phi(t)-\phi(\tau)| \bigr|\Bigr)^{-N} \notag \\ & \lesssim
  \ \lambda_j^{\frac{n+1}{2}-\frac{m}{\mu(m+2)}-\frac 2
    {m+2}-\operatorname{Re} \alpha} \bigl( \lambda_j^{-\frac{2}{m+2}}+|
  t-\tau | \bigr)^{-\frac{(m+2)(n-1)}4-\frac{m}{2\mu}}
  \notag \\ & \quad\qquad \times \Bigl( 1+ \lambda_j \big| | x-y|
  -|\phi(t)-\phi(\tau)| \big|\Bigr)^{-N}. \label{3.10}
\end{align}
Choosing $N=0$ in \eqref{3.10} gives
\begin{multline} \label{3.11}
  \|(B_j^\alpha f)(t, \tau, \cdot)\|_{L^\infty(\R^n)} \le
  \|K_j^\alpha(t, \tau, \cdot, \cdot)\|_{L^\infty_{x,y}} \|f(\tau,
  \cdot)\|_{L^1(\R^n)} \\
  \lesssim \lambda_j^{\frac{n+1}{2}-\frac{m}{\mu(m+2)}-\frac 2
    {m+2}-\operatorname{Re} \alpha} \bigl(
  \lambda_j^{-\frac{2}{m+2}}+| t-\tau |
  \bigr)^{-\frac{(m+2)(n-1)}{4}-\frac{m}{2\mu}} \|f(\tau,
  \cdot)\|_{L^1(\R^n)}.
\end{multline}
Interpolation between \eqref{3.9} and \eqref{3.11} yields \eqref{3.6}
in case $1\leq p\leq 2$ which
completes the proof of estimate \eqref{3.6}.
\end{proof}



\subsection*{Proof of Theorem~\ref{thm3.1}}

Now we return to the proof of Theorem~\ref{thm3.1}. From \eqref{3.7},
we have
\begin{equation} \label{3.12}
  W_j^\alpha f(t,x) = \int_0^t (B_j^\alpha f)(t,\tau,x)\, d\tau.
\end{equation}
Using Minkowski's inequality and  estimate \eqref{3.6},
we thus have that
\begin{multline} \label{3.13}
  \|W_j^\alpha f(t, \cdot)\|_{L^{p'}(\R^n)} \lesssim
  \lambda_j^{\left(\frac{1}{p}-\frac{1}{2}\right)(n+1)-\frac m {\mu
      (m+2)}-\frac 2 {m+2} -\operatorname{Re} \alpha} \\ \times
  \ \int_0^{\infty} \bigl( \lambda_j^{-\frac 2
    {m+2}}+|t-\tau|\bigr)^{-(m+2)
    \left(\frac{1}{p}-\frac{1}{2}\right)\frac{n-1}{2} -
    \frac{m}{2\mu}} \|f(\tau, \cdot)\|_{L^p(\R^n)}\, d\tau.
\end{multline}


\subsection*{1) \ Case \boldmath $\max\{p_1,1\} <p \le 2$.}

In this case, we have $1< r <2.$ Note that
\[
  \frac 1 r - \frac 1 {r'}=- (m+2)\left(
  \frac{1}{p}-\frac{1}{2}\right) \frac{n-1}{2} -\frac{m}{2\mu}+1.
\]
Then it follows from the Hardy-Littlewood-Sobolev theorem and
\eqref{3.13} that estimate \eqref{3.2} holds.



\subsection*{2) \ Case \boldmath $p_1>1$ and $1 <p <p_1$.}

In this case,
\[
(m+2)\left(\frac 1 p -\frac 1 2\right)\frac {n-1} 2 + \frac m {2\mu}
>1.
\]
Thus,
\[
\sup \limits_{t>0} \int_0^{\infty} \bigl( \lambda_j^{-\frac 2
  {m+2}}+|t-\tau|\bigr)^{-(m+2)
  \left(\frac{1}{p}-\frac{1}{2}\right)\frac{n-1}{2} -
  \frac{m}{2\mu}}\, d\tau < \infty,
\]
which together with Schur's lemma and \eqref{3.13} yields
\eqref{3.4}. \hfill \qed

\bigskip

We would like to stress that in the proof of Theorem~\ref{thm3.1} only
condition \eqref{2.21} on the function $b \in C^\infty(\R_+$ $\times
\R_+ \times \R^n)$ was used, whereas the conditions \eqref{2.22} and
\eqref{2.23} were not required,

\begin{remark} \label{rem3.3}
Note that the adjoint operator $(W_j^\alpha)^\ast$ of
$W_j^\alpha$ is of the form
\begin{equation}\label{3.14}
  (W_j^\alpha)^\ast f (t,x) = \int_t^\infty \int_{\R^n} e^{i \left( x
    \cdot \xi + (\phi(t)-\phi(\tau))|\xi|\right)} \ \overline{b_j(\tau, t,
  \xi)} \hat{f}(\tau, \xi) \frac{\dbar\xi}{|\xi|^\alpha}\, d\tau.
\end{equation}
By duality, we infer from Theorem~\ref{thm3.1} that
\begin{equation} \label{3.15}
  \|(W_j^\alpha )^\ast f\|_{L^{r'}_t L^{p'}_x (\R^{1+n}_+ ) } \lesssim
  \lambda_j^{\left(\frac{1}{p}-\frac{1}{2}\right)(n+1)-\frac m {\mu
      (m+2)} -\frac 2 {m+2} -\operatorname{Re} \alpha} \|f\|_{L^r_t
    L^p_x(\R^{1+n}_+ )}
\end{equation}
if $\max\{p_1, 1\} <p \le 2$
and
\begin{equation} \label{3.16}
  \| (W_j^\alpha f)^\ast \|_{L^{2}_t L^{p'}_x (\R^{1+n}_+ ) } \lesssim
  \lambda_j^{n(\frac{2}{p}-1 ) -\frac 4 {m+2} -\operatorname{Re}
    \alpha} \|f\|_{L^2_t L^p_x(\R^{1+n}_+ ) }.
\end{equation}
if  $p_1 >1$ and $1<p<p_1$. Here, $r$ is given in \eqref{3.1}.
\end{remark}

\smallskip

As an application of Theorem~\ref{thm3.1}, we obtain the boundedness
of the operator $W_j^\alpha$ from $L^r_t L^p_x$ to $L^\infty_t L^2_x$,
where $1<r, p \le 2$.

\begin{theorem}\label{thm3.4}
Let $n \ge 2$ and $\mu\ge \max\{2,m/2\}$. Then\/\textup{:}

\noindent \textup{(i)} \ For $\ds \max\{ p_1,1 \} <p\le 2$ and
$r$ be as in \eqref{3.1}, we have that
\begin{equation} \label{3.17}
  \| W_j^\alpha f \|_{L^\infty_t L^{2}_x(\R^{1+n}_+)} \lesssim
  \lambda_j^{\left(\frac{1}{p}-\frac{1}{2}\right)\frac{n+1}{2}-\frac m
    {2 \mu (m+2)} -\frac 2 {m+2}- \operatorname{Re}
    \alpha}\left\|f\right\|_{L^r_t L^p_x(\R^{1+n}_+)}.
\end{equation}
Consequently,
\begin{equation} \label{3.18}
  \| W_j^\alpha f \|_{L^\infty_t L^{2}_x( \R^{1+n}_+ )} \lesssim
  \|f\|_{L^r_t L^p_x(\R^{1+n}_+ )} \quad \textup{if }
  \operatorname{Re} \alpha=\left( \frac 1 p -\frac 1
  2\right)\frac{n+1}{2}-\frac m { 2\mu (m+2)}-\frac 2 {m+2}.
\end{equation}

\noindent
\textup{(ii)} \ For $ p_1>1$ and $1 <p <p_1$, we have that
\begin{equation}\label{3.19}
  \| W_j^\alpha f \|_{L^\infty_t L^{2}_x(\R^{1+n}_+ )} \lesssim
  \lambda_j^{n\left(\frac{1}{p}-\frac{1}{2} \right)-\frac{3}{m+2} -
    \operatorname{Re} \alpha} \|f\|_{L^2_t L^p_x(\R^{1+n}_+ )}.
\end{equation}
In particular,
\begin{equation}\label{3.20}
\| W_j^\alpha f \|_{L^\infty_t L^{2}_x(\R^{1+n}_+ )}   \lesssim  \|f\|_{L^2_t
L^p_x(\R^{1+n}_+ )} \quad \textup{if }
\operatorname{Re} \alpha=n\left(\frac{1}{p}-\frac{1}{2} \right)-\frac{3}{m+2}.
\end{equation}
\end{theorem}

\begin{proof}
For given $j \in \mathbb Z$ and $\alpha\in\mathbb C$, denote $U=W_j^\alpha
f$. Then from \eqref{2.25} we have
\begin{equation*}
  U(t) = \int_0^t e^{i\left(\phi(t)- \phi(\tau)\right) \sqrt{-\Delta}} b_j(t,
  \tau, D_x) (-\Delta)^{-\alpha/2} f(\tau)\,d\tau,
\end{equation*}
where $b_j(t, \tau, D_x)$ is the pseudodifferential operator with full
symbol $b_j(t, \tau, \xi)$. Then $U(t)$ solves the Cauchy problem
\begin{equation*}
\left\{ \enspace
\begin{aligned}
& i \p_t U(t) = -t^{m/2} \sqrt{-\Delta} U(t) + i b_j(t, t,
  D_x)(-\Delta)^{-\alpha/2} f(t) \\
& \hspace{35mm} + i \int_0^t e^{i \left(\phi(t)- \phi(\tau)\right)
    \sqrt{-\Delta}} \p_t b_j(t, \tau, D_x) (-\Delta)^{-\alpha/2}
  f(\tau)\,d\tau , \\
& U(0)=0.
\end{aligned}
\right.
\end{equation*}
Multiplying by $\overline{U(t)}$ and then integrating over $\R^n$
yields
\begin{align*}
i \left\langle\p_t U(t), U(t)\right\rangle &=- t^{m/2}
\bigl\langle\sqrt{-\Delta } U(t), U(t)\bigr\rangle + i \left\langle
b_j(t, t, D_x)(-\Delta)^{-\alpha/2} f(t), U(t) \right\rangle \\
& \hspace{20mm} + i \left\langle\int_0^t e^{i (\phi(t)- \phi(\tau))
  \sqrt{-\Delta}} \p_t b_j(t, \tau, D_x) (-\Delta)^{-\alpha/2}
f(\tau)\,d\tau, U(t)\right\rangle,
\end{align*}
and, therefore,
\begin{multline*}
  \frac 1 2 \frac d {dt} \|U(t)\|^2 = \operatorname{Re} \left\langle
  \int_0^t e^{i (\phi(t)- \phi(\tau)) \sqrt{-\Delta}} \p_t b_j(t,
  \tau, D_x) (-\Delta)^{-\alpha/2} f(\tau)\,d\tau, U(t) \right\rangle \\
  + \operatorname{Re}\left\langle b_j^\ast (t, t, D_x)
  (-\Delta)^{-\alpha/2} U(t), f(t)\right\rangle.
\end{multline*}
 Consequently,
\begin{align*}
  \|U(s)\|^2 &= 2 \operatorname{Re} \int_0^s \left\langle \int_0^t
  e^{i (\phi(t)- \phi(\tau)) \sqrt{-\Delta}} \p_t b_j(t, \tau, D_x)
  (-\Delta)^{-\alpha/2} f(\tau)\,d\tau, U(t) \right\rangle\,dt \notag
  \\
  & \quad\qquad + 2 \operatorname{Re}\int_0^s \left\langle b_j^\ast(t,
  t, D_x) (-\Delta)^{-\alpha/2} U(t), f(t) \right\rangle \,dt \notag
  \\
  & \lesssim  \left| \int_0^s \int_{\R^n} L_j^\alpha f
  (t,x) \overline{W_j^\alpha f (t,x)}\, dx dt \right| \notag \\
  & \quad\qquad + \left| \int_0^s \int_{\R^n} b_j^\ast(t, t, D_x)
  W_j^{2\alpha} f (t,x) \overline{f(t,x)}\, dx dt \right| \notag
  \\ &= \textup{I} +\textup{II}, 
\end{align*}
where
\begin{align*}
  \textup{I}&= \left| \int_0^s \int_{\R^n} L_j^\alpha f (t,x)
  \overline{ W_j^\alpha f (t,x)} \, dx dt \right| \\
  \textup{II}&= \left| \int_0^s \int_{\R^n} b_j^\ast(t, t, D_x) W_j^{2\alpha}
  f (t,x) \overline{f(t,x) }\,dx dt \right|, \\
\intertext{and}
  L^\alpha_j f(t,x) &= \int_0^t\int_{\R^n} e^{i \left( x
  \cdot \xi + (\phi(t)-\phi(\tau))|\xi|\right)} \p_t b_j(t, \tau, \xi)
\hat{f}(\tau, \xi) \,\frac{\dbar\xi}{|\xi|^\alpha} \,d\tau.
\end{align*}
From \eqref{2.21}, one has that, for any fixed $t>0$, $b_j(t, t, D_x)
\in \Psi^{-\frac 2 {m+2}}(\R^n)$, and then $b_j^\ast(t, t, D_x) \in
\Psi^{-\frac 2 {m+2}}(\R^n)$, which yields that the term $\textup{II}$
is essentially
\[
\left| \int_0^s\int_{\R^n} (W_j^{2\alpha+\frac 2{m+2}} f )(t,x)
\overline{f(t,x) }\, dx dt \right|,
\]
and thus by application of Theorem~\ref{thm3.1} it follows that
\begin{equation} \label{3.23}
\textup{II}  \lesssim \left\{
\begin{aligned}
  & \lambda_j^{(n+1) \left(\frac{1}{p}-\frac{1}{2} \right) -\frac m {\mu (m+2)}
    -\frac 4 {m+2}- 2\operatorname{Re} \alpha}\,\|f\|_{L^r_t
    L^p_x(\R^{1+n}_+ )}^2 && \text{if } \max\{p_1, 1\} < p\le 2, \\
  & \lambda_j^{n\left(\frac{2}{p}-1\right) -\frac 6 {m+2} - 2\operatorname{Re}
    \alpha}\, \|f\|_{L^2_t L^p_x( \R^{1+n}_+ )}^2 && \text{if } p_1 >1
  \text{ and } 1<p<p_1.
\end{aligned}
\right.
\end{equation}
As for the term $\textup{I}$, note that
\[
  \textup{I} = \left| \int_0^s \int_{\R^n} ( W_j^\alpha)^\ast
  L_j^\alpha f (t,x) \overline{ f (t,x)}\, dx dt \right| \le
  \left\|(W_j^\alpha)^{\ast} L_j^\alpha f \right\|_{L^{\rho'}_t L^{p'}_x(
    \R^{1+n}_+ )} \|f\|_{L^\rho_t L^p_x(\R^{1+n}_+ )}.
\]
For any $t>0$, we have from \eqref{3.14} that
\begin{multline} \label{3.24}
  (W_j^{\alpha} )^\ast  L_j^{\alpha} f (t,x) \\
\begin{aligned}
  & = \int_t^\infty \int_0^\tau \int_{\R^n} e^{i \bigl(x \cdot \xi +
    (\phi(t)-\phi(s)) |\xi|\bigr)} \overline{b_j(\tau, t, \xi)}
  \p_\tau b_j(\tau, s, \xi ) \hat{f}(s, \xi)
  \frac{\dbar\xi}{|\xi|^{2\alpha}}\,ds\, d\tau \\
  & = \int_0^t \int_{\R^n} e^{i \left(x \cdot \xi + (\phi(t)-\phi(s))
    |\xi|\right)} \Bigl( \int_t^\infty \overline{b_j(\tau, t, \xi)}
  \p_\tau b_j(\tau, s, \xi ) \, d\tau \Bigr) \hat{f}(s, \xi)
  \frac{\dbar\xi}{|\xi|^{2\alpha}}\,ds \\
  & \qquad \qquad + \int_t^\infty \int_{\R^n} e^{i \left(x \cdot \xi +
    (\phi(t)-\phi(s)) |\xi|\right)} \Bigl( \int_s^\infty
  \overline{b_j(\tau, t, \xi)} \p_\tau b_j(\tau, s, \xi )\, d\tau
  \Bigr)\hat{f}(s, \xi) \frac{\dbar\xi}{|\xi|^{2\alpha}}\,ds.
\end{aligned}
\end{multline}
Due to conditions \eqref{2.21}-\eqref{2.23}, one has that the first
and second term in \eqref{3.24} are essentially $W_j^{2\alpha + \frac
  2 {m+2} } f$ and $\bigl(W_j^{2\alpha+ \frac 2 {m+2}}\bigr)^\ast f$ ,
respectively, where $b \in C^\infty(\R_+ \times \R_+ \times \R^n)$
satisfies condition \eqref{2.21}. Then, by applying
Theorem~\ref{thm3.1} and estimates \eqref{3.15} and \eqref{3.16}, we
have that
\begin{equation*}
\textup{I}  \lesssim \left\{
\begin{aligned}
  & \lambda_j^{(n+1) \left(\frac{1}{p}-\frac{1}{2} \right) -\frac m {\mu (m+2)}
    -\frac 4 {m+2}- 2\operatorname{Re} \alpha}\|f\|_{L^r_t
    L^p_x(\R^{1+n}_+ )}^2 && \text{if } \max\{p_1, 1\} < p\le 2, \\
  & \lambda_j^{n\left(\frac{2}{p}-1 \right) -\frac 6 {m+2} - 2\operatorname{Re}
    \alpha} \|f\|_{L^2_t L^p_x( \R^{1+n}_+ )}^2 && \text{if } p_1 >1
  \text{ and } 1<p<p_1,
\end{aligned}
\right.
\end{equation*}
which together with \eqref{3.23} yields that
\begin{align*}
  \|U(t)\|^2 \lesssim \left\{
\begin{aligned}
  & \lambda_j^{(n+1) (\frac{1}{p}-\frac{1}{2} ) -\frac m {\mu (m+2)}
    -\frac 4 {m+2}- 2\operatorname{Re} \alpha}\|f\|_{L^r_t
    L^p_x(\R^{1+n}_+ )}^2 && \text{if } \max\{p_1, 1\} < p\le 2, \\
  & \lambda_j^{n(\frac{2}{p}-1 ) -\frac 6 {m+2} - 2\operatorname{Re}
    \alpha} \|f\|_{L^2_t L^p_x(\R^{1+n}_+ )}^2 && \text{if } p_1 >1
  \text{ and } 1<p<p_1.
\end{aligned}
\right.
\end{align*}
Note that $\|W_j^\alpha f (t, \cdot)\|_{L^2(\R^n)} =
\|U(t)\|$. Therefore, we have obtained estimates
\eqref{3.17}-\eqref{3.20} which completes the proof of
Theorem~\ref{thm3.4}.
\end{proof}

\begin{remark}  \label{rem3.5}
With similar arguments as in the proof of Theorem~\ref{thm3.4}, we
have from Theorem~\ref{thm3.1} and estimates \eqref{3.15} and
\eqref{3.16} that the operator $(W_j^\alpha)^\ast$ also satisfies the
estimates \eqref{3.17}-\eqref{3.20}.
\end{remark}

\smallskip

Note that if $r=p$ for $r$ defined in \eqref{3.1}, then $r=p=p_0$.
Combining Theorem~\ref{thm3.1} and the kernel estimate \eqref{3.10},
we obtain boundedness of the operator $W_j^\alpha$ from
$L^{p_0}(\R^{1+n}_+)$ to $L^q(\R^{1+n}_+)$ for certain
$\alpha
\in\mathbb C$ when $q_0 \le q \le \infty$.

\begin{theorem}\label{thm3.6}
Let $\mu \ge \max\{2,m/2\}$ and $q_0 \le q \le \infty$. Then
\begin{equation} \label{3.25}
  \|W_j^\alpha f\|_{L^q(\R^{1+n}_+)}\lesssim \|f\|_{L^{p_0}(\R^{1+n}_+)},
\end{equation}
where
\begin{equation*}
  \operatorname{Re}\alpha=n-\frac 2 {m+2} -\left( n+\frac 2
               {m+2}\right) \left( \frac 1 q +\frac 1{q_0} \right).
\end{equation*}
%
\end{theorem}

\begin{proof}
(i) \ \textit{Case  $q=q_0$.} \  Note that
\[
 n
  -\frac 2 {q_0} \left( n+\frac 2 {m+2} \right)=  \left(\frac 1 {p_0} -\frac 1 2 \right)(n+1) -\frac m {\mu (m+2)}. \]
An application of \eqref{3.3} with $r=p$ yields that
\begin{equation} \label{3.27}
  \|W_j^\alpha f\|_{L^{q_0}(\R^{1+n}_+)}\lesssim
  \|f\|_{L^{p_0}(\R^{1+n}_+)}, \quad \operatorname{Re}\alpha=n-\frac 2
       {m+2} -\frac 2 {q_0} \left( n+\frac 2 {m+2} \right).
\end{equation}
\smallskip

\noindent
(ii) \ \textit{Case $q=\infty$.} \ In order to derive \eqref{3.25}, it
suffices to show that the integral kernel $K_j^\alpha$ defined in
\eqref{3.7} satisfies
\begin{equation} \label{3.28}
\sup_{(t,x)\in \R_+^{1+n}} \int_{\R^{1+n}_+}|K_j^\alpha(t, \tau, x,
y)|^{q_0} d\tau dy <\infty, \quad \operatorname{Re}\alpha=n-\frac 2
{m+2} -\frac 1 {q_0} \Bigl( n+\frac 2 {m+2} \Bigr).
\end{equation}
In fact, from \eqref{3.7} we have
\[
W_j^\alpha f  (t,x) = \int_0^t \int_{\R^{n}} K_j^\alpha(t, \tau, x,
y) f(\tau, y)\, dy d\tau.
\]
By H\"older's inequality,   then
\begin{equation} \label{3.29}
  \|W_j^\alpha f\|_{L^\infty(\R^{1+n}_+)} \lesssim
  \|f\|_{L^{p_0}(\R^{1+n}_+)}, \quad \operatorname{Re}\alpha=n-\frac 2
       {m+2} -\frac 1 {q_0} \Bigl( n+\frac 2 {m+2} \Bigr).
\end{equation}
Now it remains to derive estimate \eqref{3.28}. In fact, due to the
kernel estimate \eqref{3.10}, for any $N>n$ and $\alpha \in \mathbb C$
with $\operatorname{Re}\alpha=n-\frac 2 {m+2} -\frac 1 {q_0} \bigl(
n+\frac 2 {m+2} \bigr)$, we have by \eqref{3.10}
\begin{multline*}
\int_{\R^{1+n}_+}  |K_j^\alpha(t, \tau, x, y)|^{q_0} \,d\tau dy \\
\begin{aligned}
& \lesssim\lambda_j^{\left(\frac{n+1}{2}-\operatorname{Re}\alpha -\frac m
    {\mu (m+2)}-\frac 2 {m+2} \right)q_0} \int_0^\infty \bigl(
  \lambda_j^{-\frac{2}{m+2}}+| t-\tau | \bigr)^{-
    \left(\frac{(m+2)(n-1)}{4}+\frac{m}{2\mu}\right)q_0}\, d\tau\\
& \qquad \times \int_{\R^n} \bigl( 1+ \lambda_j\big||x-y|
  -|\phi(t)-\phi(\tau)| \big|\bigr)^{-N}\, dy \\
& \lesssim \lambda_j^{\left(\frac{n+1}{2}-\operatorname{Re}\alpha -\frac m
    {\mu (m+2)}-\frac 2 {m+2} \right)q_0} \int_0^\infty \bigl(
  \lambda_j^{-\frac{2}{m+2}}+| t-\tau |
  \bigr)^{-\left(\frac{(m+2)(n-1)}{4}+\frac{m}{2\mu}\right)q_0} \,
  d\tau\\
& \qquad \times \lambda_j^{-1} \int_0^\infty (1+ r)^{-N} \bigl(
  \lambda_j^{-1} r+ |\phi(t)-\phi(\tau)|\bigr)^{n-1}\, dr \\
& = \lambda_j^{\left(\frac{n+1}{2}-\operatorname{Re}\alpha -\frac m {\mu
      (m+2)} -\frac 2 {m+2}\right) q_0 -1} \\
& \qquad \times \int_0^\infty \bigl( \lambda_j^{-\frac{2}{m+2}}+|
  t-\tau | \bigr)^{-
    \left(\frac{(m+2)(n-1)}{4}+\frac{m}{2\mu}\right)q_0} \bigl(
  \lambda_j^{-1}+ |\phi(t)-\phi(\tau)|\bigr)^{n-1} \,d\tau\\
& \qquad \times \int_0^\infty ( 1+ r)^{-N} \left(\frac{ r+ \lambda_j
    |\phi(t)-\phi(\tau)|}{1+\lambda_j|\phi(t)-\phi(\tau)|}
  \right)^{n-1}\, dr \\
& \lesssim \lambda_j^{\left(\frac{n+1}{2}-\operatorname{Re}\alpha -\frac m
    {\mu (m+2)}-\frac 2 {m+2} \right) q_0 -1} \int_0^\infty \bigl(
  \lambda_j^{-\frac{2}{m+2}}+| t-\tau |
  \bigr)^{-\left(\frac{(m+2)(n-1)}{4}+\frac{m}{2\mu}\right)q_0+\frac{(m+2)(n-1)}{2}}\, d\tau\\
& \lesssim \lambda_j^{\left(n-\operatorname{Re}\alpha -\frac 2
    {m+2} \right)q_0-n-\frac{2}{m+2}}=1,
\end{aligned}
\end{multline*}
and hence \eqref{3.28} holds.

\smallskip
\noindent
(iii) \ \textit{Case $q_0<q <\infty$.} \ Applying Stein's
interpolation theorem, one obtains that estimate \eqref{3.25} holds by
interpolating between estimates \eqref{3.27} and \eqref{3.29}.
\end{proof}

\smallskip

Now we consider boundedness of the operator $W_j$ from $L^r_t
L^p_x(S_T)$ to $L^s_t L^q_x(S_T)$, where
$1/p$ is symmetric around $1/p_0$.


\begin{theorem}\label{thm3.7}
Let $n \ge 2$. Further let $p_1<p<p_2$ if $n=2$, $m \ge 2$ or if $n
\ge 3$, and $1<p<\frac{7\mu}{4\mu -2}$ if $n=2$, $m=1$. Then, for any
$\mu \ge \mu_\ast$ and $T>0$,
\begin{equation} \label{3.30}
 \|W_j f\|_{L^s_t L^q_x(S_T)}\lesssim \|f\|_{L^r_t
   L^p_x(S_T)},
\end{equation}
where $r$ is defined as in \eqref{3.1} and
\begin{equation} \label{3.31}
\left\{ \enspace
\begin{aligned}
&  \frac 1 q= \frac 1 p - \frac 4 {(m+2)(n+1)} \left( 1+ \frac
  m {2\mu} \right),   \\
& \frac 1 s= \frac{(m+2)(n-1)}4 \left( \frac 1 2 -\frac 1 q\right)
  +\frac m {4 \mu} .
\end{aligned}
\right.
\end{equation}
\end{theorem}

\begin{proof}
Since $1/p$ is symmetric around $1/p_0$, by duality
it suffices to consider the case $\max\{p_1, 1\} <p\le p_0$.

In order to derive \eqref{3.30}, we now need a further dyadic
decomposition with respect to the time variable $t$. Choose a
function $\eta \in C^\infty_c(\R_+)$ such that $0\le \eta \le 1$,
$\operatorname{supp} \eta \subseteq [1/2, 2]$, and
\[
\sum \limits_{ \ell =-\infty}^\infty  \eta(2^{- \ell } t) =1.
\]
Let us fix $\lambda=2^j$ and set
\begin{equation*}
\eta_0(t)=\sum_{k \le 0} \eta(\lambda 2^{- k}t), \quad
\eta_\ell(t)=\eta(\lambda 2^{-\ell}t) \quad \text{for } \ell \in \N.
\end{equation*} Then,
\[
W_j   f(t, x)=\sum \limits_{k=0}^\infty G_k f (t,x),
\]
where
\begin{equation} \label{3.32}
G_k f (t,x) = \int_0^t\int_{\R^n} e^{i \left( x \cdot \xi +
(\phi(t)-\phi(\tau))|\xi|\right)} \eta_k(t-\tau)\, b_j(t, \tau, \xi)
\hat{f}(\tau, \xi) \ \dbar\xi d\tau.
\end{equation}
Hence, to derive \eqref{3.30}, it suffices to show that, for any $ k
\in \N_0$,
\begin{equation} \label{3.33}
\|G_k f\|_{L^s_t L^q_x(S_T)}\lesssim  2^{-\varepsilon_p k}
\|f\|_{L^r_t L^p_x(S_T)}
\end{equation}
for some $\varepsilon_p>0$.
From  \eqref{3.1}  and \eqref{3.31},  we know that
\begin{equation*}  
\frac{(m+2)n}{2} \left( \frac 1 p -\frac 1 q \right) + \frac 1 r- \frac
1 s =2.
\end{equation*}
Due to scaling invariance, we need to consider only the case
$\lambda=1$ (by a change of variable if $\lambda \neq 1$). Repeating
the arguments which are used to prove \eqref{3.2}, we get that, for any
$k \in \N_0$,
\begin{equation}\label{3.34}
\|G_k f \|_{L^{r'}_t L^{p'}_x (S_T)}\lesssim  2^{-k \left(
(m+2)\left(\frac 1 p -\frac 1 2\right)\frac {n-1}{2} + \frac m {2\mu}\right)}
\|f\|_{L^r_t L^p_x(S_T)}.
\end{equation}
Note that $(m+2)\left(\frac 1 p -\frac 1 2\right)\frac {n-1}{2}+\frac
m {2\mu}>\frac 1 3 $, since $p \le p_0$.

Furthermore, an immediate consequence of \eqref{3.17} for $\alpha= 0$
is
\begin{equation*}
\|G_k f \|_{L^\infty_t L^{2}_x(S_T)}\lesssim  \|f\|_{L^r_t
L^p_x(S_T)},
\end{equation*}
and thus, for any $1< \rho<\infty$,
\begin{equation} \label{3.35}
\|G_k f \|_{L^\rho_t L^{2}_x (S_T)}\lesssim \|f\|_{L^r_t L^p_x(S_T)}.
\end{equation}

Choose
\begin{equation} \label{3.36}
  \theta=\frac{4p(2\mu +m)}{\mu(m+2)(n+1)(2-p)} -1.
\end{equation}
Then $ 0\le \theta \le 1$ and, for the number $q$ from \eqref{3.31},
\[
\frac 1 q= \frac{\theta}{p'} +\frac{1-\theta}{2}.
\]
For $s$ from \eqref{3.31} and $\theta$ from \eqref{3.36}, we define
$s_0$ by
\[
2 \left( \frac 1 s -\frac 1 {s_0} \right)= \theta \left( (m+2)\left(\frac 1
p -\frac 1 2\right)\frac {n-1}{2} + \frac m {2\mu} \right)
\]
and then set $\rho=\rho_\ast$ such that
\[
\frac 1 {s_0}= \frac \theta {r'} +\frac{1-\theta}{\rho_\ast}.
\]
Since $2<s<s_0$, by interpolating between \eqref{3.34} and
\eqref{3.35} when $\rho=\rho_\ast$, we obtain that
\begin{equation} \label{3.37}
 \|G_k f\|_{L^{s_0}_t L^q_x(S_T)}\lesssim 2^{-2 k \left(\frac 1 s -\frac 1
   {s_0}\right)}\|f\|_{L^r_t L^p_x(S_T)}.
\end{equation}

Let $\{I_\ell\}$ be non-overlapping intervals of side length $2^{k}$
and $\bigcup_\ell I_\ell =\R_+$, and denote by $\chi_I$ the
characteristic function of $I$. In view of \eqref{3.32} and the
compact support of $\eta_k$, we have that if $f(t, x)= 0 $ for $t
\notin I_\ell$, then $G_k f (t,x) = 0 $ for $t \notin I^\ast_\ell$,
where $I^\ast_\ell$ is the interval with the same center as $I_\ell$
but of side length $C_0 2^{k}$ with some constant $C_0=C_0(\eta) >0$.
Thus, from Minkowski's inequality
\begin{equation} \label{3.38}
  \|G_k f (t, \cdot)\|_{L^q(\R^n)}^s \le \Bigl( \sum \limits_{\ell}
  \|G_k(\chi_{I_\ell } f ) (t, \cdot)\|_{L^q(\R^n)} \Bigr)^s \lesssim \sum
  \limits_{\ell}\|G_k(\chi_{I_\ell}f) (t, \cdot)\|_{L^q(\R^n)}^s,
\end{equation}
Denote $\overline{I_\ell^\ast}=I_\ell^\ast \cap (0, T)$. Estimate
\eqref{3.38} together with H\"older's inequality and \eqref{3.37}
yields that, for any $ k \in \N_0$,
\begin{align*}
  \|G_k f\|_{L^s_t L^q_x(S_T)}^s & \lesssim \sum
  \limits_{\ell}\|G_k(\chi_{I_\ell}f)\|_{L^s_t
    L^q_x(\overline{I_\ell^\ast} \times \R^n)}^s\\
  & \lesssim \sum
  \limits_{\ell}|\overline{I_\ell^\ast}|^{1-s/{s_0}}\|G_k(\chi_{I_\ell}f)\|_{L^{s_0}_t
    L^q_x(\overline{I_\ell^\ast} \times \R^n)}^s \\
  & \lesssim 2^{k \left(1 -s/{s_0}\right)}2^{-2ks \left(1/s -1/{s_0}\right)} \sum
  \limits_{\ell} \|\chi_{I_\ell} f\|_{L^r_t L^p_x(S_T)}^s \\
  & \lesssim 2^{-k \left(1 -s/{s_0}\right)} \|f\|_{L^r_t L^p_x(S_T)}.
\end{align*}
Therefore, we get estimate \eqref{3.33} with $\varepsilon_p=1 - s /{s_0}$
and, hence, \eqref{3.30} holds.
\end{proof}


By a similar argument as in the proof of Theorem~\ref{thm3.7}, we
obtain the boundedness of operator $W_j$ from $L^2_t L^p_x(S_T)$ to
$L^s_t L^q_x(S_T)$ when $p_1>1$ and $1<p<p_1$.

\begin{theorem}\label{thm3.8}
Let $n\geq3$ or $n =2$, $m\geq2$. Suppose $1<p<p_1$. Then, for $\mu
\ge \max\{2,mn/2\}$ and $T>0$, we have that
\begin{equation}  \label{3.39}
 \|W_j f\|_{L^s_t L^q_x(S_T)}\lesssim \|f\|_{L^2_t L^p_x(S_T )},
\end{equation}
where
\begin{equation} \label{3.40}
\left\{ \enspace
\begin{aligned}
  & \frac 1q=\frac {2n}{p\left(n+1\right)} -\frac{n-1}{2(n+1)}-\frac
  {m+6\mu} {\mu(m+2)(n+1)}, \\
  & \frac 1 s= (m+2) \left( \frac 1 2 -\frac 1 q\right)
  \left(\frac{n-1}4 \right)+\frac m {4 \mu}.
\end{aligned}
\right.
\end{equation}
\end{theorem}

\begin{proof}
Note that when $1< p<p_1$, we have
\[
  (m+2)\left(\frac 1 p -\frac 1 2\right)\left(\frac
  {n-1}{2}\right)+\frac m {2\mu}>1.
\]
Then we can apply similar arguments as in the proof
Theorem~\ref{thm3.7} to obtain \eqref{3.39}. We omit the details.
\end{proof}

\begin{remark} \label{rem3.9}
By similar arguments as above one can show that adjoints $(W_j)^\ast$
of $W_j$ also satisfy estimates \eqref{3.30} and \eqref{3.39},
respectively, under assumptions \eqref{3.31} and \eqref{3.40}.
\end{remark}

\section{Mixed-norm estimates for the linear generalized Tricomi equation}\label{sec4}

In this section, based on the mixed-norm space-time estimates of the
Fourier integral operators $W_j^\alpha$ obtained in
Section~\ref{sec3}, we shall establish Strichartz-type estimates for
the linear generalized Tricomi equation.

\smallskip

First we consider the inhomogeneous equation, i.e., problem
\eqref{2.3}.

\begin{theorem}\label{thm4.1}
Let $n \ge 2$. Suppose $w$ is a solution of \eqref{2.3} in $S_T$ for
some $T>0$. Then\/\textup{:}

\smallskip
\noindent
\textup{(i)} \ For $\mu \ge \mu_\ast$,
\begin{equation}\label{4.1}
 \|w\|_{L^s_t L^q_x(S_T)}\lesssim \|f\|_{L^r_t L^p_x(S_T)}
\end{equation}
provided that $p_1 <p<p_2$ if $n\geq3$ or $n=2$, $m\geq2$, and
$1<p<7\mu/\!\left(4\mu -2\right)$ if $n=2$ and $m=1$. Here  $r=r(p,\mu)$ is
as in \eqref{3.1} and $q$ and $s$ are taken from \eqref{3.31}.

\smallskip

\noindent
\textup{(ii)} \ For $\mu \ge \max\{2,m/2\}$,
\begin{equation} \label{4.2}
\|w\|_{L^q(S_T)}\lesssim \left\| |D_x|^{\gamma-\gamma_0}
f\right\|_{L^{p_0}(S_T)}, \quad q_0 \le q <\infty,
\end{equation}
where
\begin{equation} \label{4.3}
\left\{ \enspace
\begin{aligned}
&\gamma=\gamma(m,n,q)= \frac n 2 -\frac 1 q \left( n +\frac 2
  {m+2}\right),\\
&\gamma_0=\gamma_0(m,n,\mu)=\frac 1 {q_0} \left( n + \frac 2 {m+2}
  \right) +\frac 2 {m+2} -\frac n 2.
\end{aligned}
\right.
\end{equation}

%

\smallskip

\noindent
\textup{(iii)} \ For $\mu \ge \max\{2,m/2\}$, $\max\{p_1, 1\}< p \le
2$, and $0\leq t \leq T$,
\begin{equation} \label{4.4}
\|w( t,
\cdot)\|_{\dot{H}^\gamma(\R^n)}  \lesssim  \| f\|_{L^r_t L^p_x(S_T)},
\end{equation}
where $r=r(p,\mu)$ is defined in \eqref{3.1} and
\[
\gamma=\gamma(m,n,\mu,p)=\frac 2 {m+2}+\frac m {2\mu (m+2)} - \left(
  \frac 1 p -\frac 1 2 \right) \frac {n+1} 2.
\]

\smallskip

\noindent
\textup{(iv)} \ For $\mu \ge \max\{2,m/2\}$, $\gamma \in \R$, and
$0\leq t\leq T$,
\begin{equation} \label{4.5}
 \|w( t, \cdot)\|_{\dot{H}^\gamma(\R^n)} \lesssim \left\|
 |D_x|^{\gamma-\gamma_0} f\right\|_{L^{p_0}(S_T)},
\end{equation}
where $\gamma_0$ is from \eqref{4.3}.
\end{theorem}

\smallskip

\begin{remark} \label{rem4.2}
If we choose $\mu=\mu_\ast$, then
\[
p_0=p_0^\ast= \frac{2\mu_*}{\mu_*+2},
\quad
q_0=q_0^\ast=\frac{2\mu_*}{\mu_*-2},
\]
and for $\gamma$ and $\gamma_0$ defined in \eqref{4.3},
\[
  \gamma(m,n,q_0^*)=\gamma_0(m,n,\mu_*)=\frac 1{m+2}.
\]
Thus, we have from \eqref{4.2} that
\[
\|w\|_{L^{q_0^\ast}(S_T)}\lesssim  \|  f\|_{L^{p_0^\ast}(S_T)},
\]
which, for any $\rho\in\R$, together with $\bigl[|D_x|^\rho, \p_t^2 -t^m
  \Delta\bigr]=0$ implies that
\[
  \left\| |D_x|^\rho w\right\|_{L^{q_0^\ast}(S_T)}\lesssim \|
  |D_x|^\rho f\|_{L^{p_0^\ast}(S_T)}.
\]
\end{remark}



\subsection*{Proof of Theorem~\ref{thm4.1}}

${ }$

\smallskip
\noindent
(i) \
One obtains \eqref{4.1} by applying Proposition~\ref{prop2.1} and Theorem~\ref{thm3.7}
directly.

\smallskip
\noindent
(ii) \ For $\alpha \in \mathbb C$, the Fourier transform of
$|D_x|^{\alpha} f(t,x)$ with respect to the variable $x$ is
$|\xi|^{\alpha} \hat{f}(t, \xi)$. Thus, we can write $W_j f$ as
\begin{equation*}
  W_j f(t,x) = \int_0^t \int_{\R^n} e^{i \left(x \cdot \xi +
    (\phi(t)-\phi(\tau)) |\xi|\right)} \Theta(|\xi|/2^j)\, b(t, \tau,
  \xi) (\widehat{|D_x|^{\alpha} f})(\tau, \xi) |\xi|^{-\alpha}
  \ \dbar\xi d\tau
\end{equation*}
and $W_j (f) = W_j^{\alpha}(|D_x|^{\alpha}f)$.
Therefore, applying Theorem~\ref{thm3.6}, we get that
\begin{equation*}
  \|W_j f\|_{L^q(S_T)}= \|W_j^{\gamma-\gamma_0} (
  |D_x|^{\gamma-\gamma_0} f)\|_{L^q(S_T)} \lesssim \|
  |D_x|^{\gamma-\gamma_0} f\|_{L^{p_0}(S_T)},
\end{equation*}
which together with Proposition~\ref{prop2.1}  yields  \eqref{4.2}.

\smallskip
\noindent
(iii) \ Note that $\left[|D_x|^\gamma, \p_t^2 -t^m \Delta\right]=0$ and then
\begin{equation}\label{4.7}
  (\p_t^2 -t^m \Delta) (|D_x|^\gamma w)=|D_x|^\gamma f.
\end{equation}
From (ii) we know that $W_j(|D_x|^\gamma f)=W_j^{-\gamma}(f)$. Thus,
for $\gamma=\frac 2 {m+2}+\frac m {2\mu (m+2)} - \left( \frac 1 p
-\frac 1 2 \right) \frac{n+1}2 $, we have from estimate \eqref{3.18}
that
\[
  \| W_j(|D_x|^\gamma f)( t, \cdot )\|_{L^2(\R^n)} = \| W_j^{-\gamma}
  f (t, \cdot)\|_{L^2(\R^n)} \lesssim \|f\|_{L^r_t L^p_x}.
\]
Thus, by \eqref{4.7} and Proposition~\ref{prop2.1} it follows that
\[
\|(|D_x|^\gamma w)(t, \cdot)\|_{L^2(\R^n)} \lesssim \|f\|_{L^r_t
  L^p_x},
\]
which together with  Plancherel's theorem implies that
\begin{equation*} 
\|w( t, \cdot)\|_{\dot{H}^{\gamma}(\R^n)}=\left\| |\xi|^\gamma \hat{w}( t,
\xi)\right\|_{L^2_{\xi}(\R^n)}= \left\| (|D_x|^\gamma w) ( t, \cdot)
\right\|_{L^2_x(\R^n)} \lesssim \|f\|_{L^r_t L^p_x},
\end{equation*}
and estimate \eqref{4.4} holds.

\smallskip
\noindent
(iv) \ From (ii) we also know that
\[
W_j (g) = W_j^{-\gamma_0}
(|D_x|^{-\gamma_0}g).
\]
In \eqref{3.1}, we have $r=p=p_0$ when $r=p$. Estimate \eqref{3.18}
for \[ \alpha=- \gamma_0=\left(\frac 1 {p_0} -\frac 1 2 \right)\frac {n+1} 2  -\frac m {2\mu (m+2)} -\frac 2 {m+2} \]
with $p=p_0$ yields that
\begin{equation*}
  \| W_j(g)( t, \cdot )\|_{L^2(\R^n)} =\left\| W_j^{-\gamma_0}
  (|D_x|^{-\gamma_0}g) ( t, \cdot )\right\|_{L^2(\R^n)} \lesssim \left\|
  |D_x|^{-\gamma_0} g\right\|_{L^{p_0}(S_T)},
\end{equation*}
and then, for $g=|D_x|^\gamma f$, where $\gamma \in \R$,
\begin{equation}  \label{4.8}
\left\| W_j(|D_x|^\gamma   f)(  t, \cdot )\right\|_{L^2(\R^n)}\lesssim
\left\| |D_x|^{\gamma-\gamma_0} f\right\|_{L^{p_0}(S_T)}.
\end{equation}
Therefore, one has from Plancherel's theorem,
Proposition~\ref{prop2.1}, \eqref{4.7}, and \eqref{4.8} that
\[
\|w( t, \cdot)\|_{\dot{H}^{\gamma}(\R^n)}=\left\| (|D_x|^\gamma w)(t,
\cdot)\right\|_{L^2(\R^n)} \lesssim \left\| |D_x|^{\gamma-\gamma_0}
f\right\|_{L^{p_0}(S_T)}
\]
Hence, estimate \eqref{4.5} holds. \qed

\bigskip

In case $n \ge 2$ and $m \ge 2$ if $n=2$, we have a more complete set
of inequalities for the solution of the linear generalized Tricomi
equation.

\begin{theorem}\label{thm4.3}
Let $n \ge 3$ or $n=2$,  $m \ge 2$. Suppose $w$ solves \eqref{2.3}
in $S_T$. Then\/\textup{:}

\noindent
\textup{(i)} \ For $ \mu \ge \max\{ 2,mn/2\}$ and $\frac 1
       {p_1} < \frac 1 p \le \frac 1 2 + \frac {m+6\mu} {2\mu
         n(m+2)}$,
\begin{equation} \label{4.9}
\|w\|_{L^s_t L^q_x(S_T)}\lesssim
\|f\|_{L^2_t L^p_x(S_T)}, \end{equation}
where $q$ and $s$
are defined in \eqref{3.40}.

\smallskip
\noindent
\textup{(ii)} \ For $ \mu \ge \max\{ 2,mn/2\} $ and $\frac 1 2 \le
\frac 1 p < \frac 1 2 + \frac {2\mu(n-3) +m(3n-1)} {\mu
  (m+2)(n^2-1)}$,
%
\begin{equation}\label{4.10}
\|w\|_{L^2_t L^{q}_x(S_T)}\lesssim   \|f\|_{L^{r}_t L^{p}_x(S_T)},
\end{equation}
where $r$ is defined in \eqref{3.1} and
\begin{equation} \label{4.11}
 \frac 1 q =\frac {n+1}{2np} + \frac {n-1}{4n} -\frac {m+6\mu}{2\mu
   (m+2)n}.
\end{equation}

\smallskip
\noindent
\textup{(iii)} \ For $ \mu \ge \max\{ 2, \frac{m}{2}\}$ and $ \ds 1 <p
<p_1$ and $ \gamma=\frac{3}{m+2}-n\bigl(\frac{1}{p}-\frac{1}{2}
\bigr)$,
\begin{equation}  \label{4.12}
\|w(  t, \cdot)\|_{\dot{H}^\gamma(\R^n)}\lesssim  \|f\|_{L^2_t
L^p_x(S_T)}.
\end{equation}
\end{theorem}

\begin{proof}
\noindent
(i) \ Note that, under these assumptions,
\[
  1<\frac{2\mu n(m+2)}{\mu n(m+2)+6\mu +m}\le p<p_1, \quad 2 \le q
  <\infty, \quad 2 \le s <\infty.
\]
Thus, we get estimate \eqref{4.9}  by applying
Proposition~\ref{prop2.1} and Theorem~\ref{thm3.8}.

\smallskip
\noindent
(ii) \ This will follow from the dual version of
Theorem~\ref{thm3.8}. Indeed, when
\[
\frac 1 2 \le \frac 1 p < \frac 1 2 + \frac {2\mu(n-3) +m(3n-1)} {\mu
  (m+2)(n^2-1)},
\]
then, for $q$ defined in \eqref{4.11},
\[
1<\frac{2\mu (m+2)n}{\mu (m+2)n+6\mu +m} \le q'<p_1
\]
and
\[
\frac 1 {p'}=\frac {2n}{q'(n+1)} -\frac{n-1}{2(n+1)}
-\frac{m+6\mu}{\mu(m+2)(n+1)},
\]
For $r$ defined by \eqref{3.1}, the conjugate exponent $r'$ can be
expressed by
\[
r'=\frac{8\mu p'}{\mu(m+2)(n-1)(p'-2)+2mp'}.
\]
Thus, from Remark~\ref{rem3.9}, we have that
\[
\| W_j^\ast f \|_{L^{r'}_t L^{p'}_x(S_T)}\lesssim \|f\|_{L^{2}_t
  L^{q'}_x(S_T)},
\]
and then, by duality,
\[
\| W_j f\|_{L^{2}_t L^{q}_x(S_T)}\lesssim \|f\|_{L^r_t
  L^p_x(S_T)}.
\]
Therefore, from Proposition~\ref{prop2.1} we have that estimate
\eqref{4.10} holds.

\smallskip
\noindent
(iii) \ Note again that $W_j(|D_x|^\gamma f)=W_j^{-\gamma}(f)$. Then,
in view of \eqref{4.7} and estimate \eqref{3.20} for $\alpha =
-\gamma=n\bigl(\frac{1}{p}-\frac{1}{2} \bigr)-\frac{3}{m+2}$, one has
that estimate \eqref{4.12} holds.
\end{proof}

Now we consider the homogeneous equation, i.e., problem \eqref{2.2}.

\begin{theorem} \label{thm4.4}
Let $n \ge 2$ and $\mu \ge \max\{2,m/2\}$. Suppose $v$ solves the
Cauchy problem \eqref{2.2}. Then\/\textup{:}

\noindent
\textup{(i)} \ For $q_0\le q <\infty$,
\begin{equation} \label{4.13}
\|v\|_{L^q(\R^{1+n}_+)} \lesssim
\|\varphi\|_{\dot{H}^{\gamma}(\R^n)} +
\|\psi\|_{\dot{H}^{\gamma-\frac 2 {m+2}}(\R^n)},
\end{equation}
where $\gamma=\frac{n}2 -
\frac{\left(m+2\right)n+2}{q\left(m+2\right)}$.

\noindent
\textup{(ii)} \ For $2 \le q <\infty$ when $n=2$ and $m=1$, and $2\le q
<q_1$ when $n \ge 2$ and $m \ge 2$ if $n=2$,
\begin{equation}\label{4.14}
  \|v\|_{L^s_t L^q_x(\R^{1+n}_+)}\lesssim
  \|\varphi\|_{\dot{H}^{\gamma}(\R^n)} +
  \|\psi\|_{\dot{H}^{\gamma-\frac 2 {m+2}}(\R^n)},
\end{equation}
where
\[
  \frac 1 s= \frac{(m+2)(n-1)}4 \left( \frac 1 2 -\frac 1 q\right)
  +\frac m {4 \mu}, \quad \gamma=\frac{n+1}{2}\left( \frac 1 2 -\frac
  1 q \right)-\frac m {2\mu(m+2)}.
\]

\noindent
\textup{(iii)} \ For $q_1 < q <\infty$ as well as $n \ge 2$ and $m \ge
2$ if $n=2$,
\begin{equation} \label{4.15}
\|v\|_{L^2_t L^q_x(\R^{1+n}_+)}\lesssim
\|\varphi\|_{\dot{H}^{\gamma}(\R^n)}  +
\|\psi\|_{\dot{H}^{\gamma-\frac 2 {m+2}}(\R^n)}, \quad
\end{equation}
where $\gamma=n\bigl( \frac 1 2 -\frac 1 q \bigr)-\frac 1 {m+2}.$
\end{theorem}

\begin{proof}
The goal is to prove that
\begin{equation} \label{4.16}
\|v\|_{L^\sigma_t L^\rho_x(\R^{1+n}_+)}\lesssim
\|\varphi\|_{\dot{H}^{\gamma}(\R^n)} + \|\psi\|_{\dot{H}^{\gamma-\frac
    2 {m+2}}(\R^n)}
\end{equation}
for certain $2 \le \sigma \le \infty$ and $ 2 \le \rho <\infty$.

Note that
\[
  t \left( 1+ \phi(t) |\xi|\right)^{-\frac{m+4}{2(m+2)}} \le \left( 1+
  \phi(t) |\xi|\right)^{-\frac{m}{2(m+2)}} |\xi|^{-\frac 2 {m+2}} \le
  \left( 1+ \phi(t) |\xi|\right)^{-\frac{m}{\mu(m+2)}} |\xi|^{-\frac 2
    {m+2}}.
\]
In order to establish \eqref{4.16}, from the expression of the
function $v$ in \eqref{2.4} together with \eqref{2.11} and \eqref{2.12} and
the estimates of $b_\ell(t, \xi) (1\le \ell \le 4)$ in \eqref{2.15}
and \eqref{2.16}, it suffices to show that
\begin{equation} \label{4.18}
  \|P \varphi\|_{L^\sigma_t L^\rho_x(\R^{1+n}_+)}\lesssim
  \|\varphi\|_{\dot{H}^\gamma(\R^n)},
\end{equation}
where the operator $P$ is of the form
\begin{equation*}
  (P \varphi)(t, x)= \int_{\R^n} e^{i \left(x \cdot \xi + \phi(t)
    |\xi|\right)} a(t, \xi) \hat{\varphi}(\xi)\ \dbar\xi
\end{equation*}
with $ a \in C^\infty(\R_+ \times \R^n)$ and, for any $(t,\xi) \in
\R_+^{1+n}$,
\begin{equation} \label{4.17}
\bigl|\partial_\xi^\beta a(t, \xi)\bigr|\lesssim  \left( 1+
\phi(t) |\xi|\right)^{-\frac{m}{\mu(m+2)}} |\xi|^{-|\beta|}.
\end{equation}
Note that $P\varphi$ can be written as
\begin{equation*}
  (P \varphi)(t, x) = \int_{\R^n} e^{i \left(x \cdot \xi + \phi(t)
    |\xi|\right)} a(t, \xi) \widehat{ |D_x|^\gamma \varphi } (\xi)
  \ \frac{\dbar\xi}{|\xi|^\gamma},
\end{equation*}
and, for $h=|D_x|^\gamma \varphi$, by Plancherel's theorem,
\[
  \|h\|_{L^2(\R^n)}= \left\| |\xi|^\gamma \hat{\varphi}
  \right\|_{L^2(\R^n)}=\|\varphi\|_{\dot{H}^\gamma( \R^n)}.
\]

Therefore, in order to prove \eqref{4.18}, it suffices to show that
the operator $T$, where
\begin{equation} \label{4.19}
  (T h) (t, x)= \int_{\R^n} e^{i \left(x \cdot \xi + \phi(t) |\xi|\right)} a(t,
  \xi) \hat{h}(\xi) \, \frac {\dbar\xi}{|\xi|^{\gamma}},
\end{equation}
is bounded from $L^2(\R^n)$ to $L^\sigma_t L^\rho_x(\R^{1+n}_+)$. By
duality, it suffices to show that the adjoint $T^\ast$ of $T$,
\begin{equation}  \label{4.20}
  (T^\ast f)(x) = \int_0^\infty \int_{\R^n} e^{i \left(x\cdot \xi -
    \phi(\tau) |\xi|\right)} \overline{a(\tau, \xi)} |\xi|^{-\gamma} \hat{f}(\tau,
  \xi) \ \dbar\xi d\tau,
\end{equation}
satisfies
\begin{equation} \label{4.21}
  \|T^\ast f \|_{L^2(\R^n)}\lesssim \|f\|_{L^{\sigma'}_t
    L^{\rho'}_x(\R^{1+n}_+)}.
\end{equation}

Note that
\begin{align*}
  \| T^\ast f\|_{L^2(\R^n)}^2 &= \int_{\R^n} (T^\ast f)(x)
  \overline{(T^\ast f)(x)}\, dx \\
  & = \int_{\R^{1+n}_+} T T^\ast f(t,x) \overline{f(t,x)}\, dt dx \le
  \|T T^\ast f\|_{L^\sigma_t L^\rho_x} \|f\|_{L^{\sigma'}_t
    L^{\rho'}_x}.
\end{align*}
Thus, in order to get \eqref{4.21}, we only need to show that
\begin{equation}  \label{4.22}
  \|T T^\ast f\|_{L^\sigma_t L^\rho_x}\lesssim \|f\|_{L^{\sigma'}_t
    L^{\rho'}_x}.
\end{equation}
From \eqref{4.19}  and  \eqref{4.20}, we have that
\begin{equation*}
  T T^\ast f (t,x) = \int_0^\infty \int_{\R^n} e^{i\left(
    x\cdot \xi+ (\phi(t)-\phi(\tau) ) |\xi|\right)} a(t, \xi)\overline{a(\tau,
  \xi)} \hat{f}(\tau, \xi) \, \frac{\dbar\xi}{|\xi|^{2\gamma}} d\tau.
\end{equation*}
By \eqref{4.17}, we further have that
\[ 
  \left|\p_{\xi}^\beta \bigl(a(t, \xi) \overline{a(\tau,\xi)}\bigr)
  \right|\lesssim \bigl( 1+ |\phi(t) - \phi(\tau)| |\xi|
  \bigr)^{-\frac{m}{\mu(m+2)}} |\xi|^{-|\beta|}.
\]
Thus, by Proposition~\ref{prop2.1}, in order to get
\eqref{4.22}, it suffices to show that
\begin{equation*}
  \|G_j f \|_{L^\sigma_t L^\rho_x}\lesssim \|f\|_{L^{\sigma'}_t
    L^{\rho'}_x},
\end{equation*}
where the operator $G_j$ is defined as
\begin{equation*} 
 G_j f (t,x)  = \int_0^\infty
\int_{\R^n}e^{i\left( x \cdot \xi+ ( \phi(t)-\phi(\tau) )|\xi|\right)}
 \Theta(|\xi|/2^j) a(t, \xi)\overline{a(\tau, \xi)}
 \hat{ f}(\tau, \xi) \, \frac{\dbar\xi}{|\xi|^{2\gamma}} d\tau.
\end{equation*}
Note that $ G_j f$ is essentially $W_j^{2\gamma- \frac 2 {m+2}}
f$. Therefore, in order to get \eqref{4.15}, it suffices to show that
\begin{equation} \label{2.4}
  \bigl\|W_j^{2\gamma - \frac 2 {m+2} } f\bigr\|_{L^\sigma_t L^\rho_x} \lesssim
  \|f\|_{L^{\sigma'}_t L^{\rho'}_x}.
\end{equation}

\smallskip

We first show \eqref{4.13}: \ For $\gamma=\frac n 2
-\frac{n(m+2)+2}{q(m+2)}$ and $q=q_0$, we have that
\begin{equation*}
  \left( 2\gamma- \frac 2 {m+2} \right)=\left(\frac 1 {p_0} -\frac 1 2
  \right) (n+1) -\frac m {\mu (m+2)} -\frac 2 {m+2}.
\end{equation*}
Thus, we have from estimate \eqref{3.3} when $r=p=p_0$ that
\begin{equation}  \label{4.26}
  \bigl\|W_j^{2\gamma - \frac 2 {m+2} }
  \bigr\|_{L^{q_0}(\R^{1+n}_+)}\lesssim \|f\|_{L^{p_0}(\R^{1+n}_+)}.
\end{equation}
On the other hand, from  \eqref{2.25} and the compact support of $\Theta$,
\begin{equation} \label{4.27}
  \bigl\|W_j^{2\gamma - \frac 2 {m+2}}
  f\bigr\|_{L^{\infty}(\R^{1+n}_+)} \lesssim \|f\|_{L^1(\R^{1+n}_+)}.
\end{equation}
By interpolation between \eqref{4.26}  and  \eqref{4.27}, we  obtain   that
\begin{equation*}  
  \bigl\| W_j^{2\gamma - \frac 2 {m+2}}
  f\bigr\|_{L^{q}(\R^{1+n}_+)}\lesssim \|f\|_{L^{q'}(\R^{1+n}_+)},
  \quad q_0 \leq q \leq \infty.
\end{equation*}
where $q'$ is the conjugate exponent $q$.  Therefore, we get estimate
\eqref{4.13}.

\smallskip

Next we derive \eqref{4.14}: \ Since
\[
  \frac 1s= \frac{(m+2)(n-1)}4 \left( \frac 1 2 -\frac 1 q\right) +\frac m {4 \mu},
\]
we can write
\[
  \frac 1 {s'}=1- \frac{(m+2)(n-1)}4 \left( \frac 1 {q'} - \frac 1 2
  \right) - \frac m {4 \mu}
\]
Thus, when $\gamma=\frac{n+1}{2}\bigl( \frac 1 2 -\frac 1 q
\bigr)-\frac m {2\mu(m+2)}$, applying estimate \eqref{3.3} for
$\max\{p_1, 1\} < q'\le 2$, we have
\begin{equation*}  
\bigl\|  W_j^{2\gamma- \frac 2 {m+2}} f\bigr\|_{L^s_t L^q_x(\R^{1+n}_+)}\lesssim
\|f\|_{L^{s'}_t L^{q'}_x(\R^{1+n}_+)},
\end{equation*}
and, therefore, estimate \eqref{4.14} holds.

\smallskip

Finally we prove \eqref{4.15}: \ When $\gamma=n\bigl( \frac 1 2 -\frac
1 q \bigr) -\frac 1 {m+2}$, we have from \eqref{3.5} that, for $p_1
>1$ and $1 <q' < p_1$,
\begin{equation*}
  \bigl\| W_j^{2\gamma- \frac 2 {m+2}} f\bigr\|_{L^2_t
    L^q_x(\R^{1+n}_+)}\lesssim \|f\|_{L^{2}_t L^{q'}_x(\R^{1+n}_+)}.
\end{equation*}
Thus, estimate \eqref{4.15}  holds.
\end{proof}

\smallskip


\smallskip

Combining Theorems~\ref{thm4.1},~\ref{thm4.3}, and~\ref{thm4.4}, we obtain the
following results:

\begin{theorem} \label{thm4.5}
Let $u$ solve the Cauchy problem \eqref{2.1} in the
strip $S_T$. Then
\begin{equation}  \label{4.30}
\|u\|_{C_t^0\dot{H}_x^\gamma(S_T)} + \|u\|_{L^s_t L^q_x ( S_T)}  \lesssim
\|\varphi\|_{\dot{H}^\gamma(\R^n)}+
\|\psi\|_{\dot{H}^{\gamma-\frac 2 {m+2}}(\R^n)} + \|f\|_{L^r_t L^p_x (S_T)}
\end{equation}
provided that the exponents $p,\, q,\, r$, and $s$ satisfy scaling
invariance condition \eqref{SI}
and one of the following sets of conditions\/\textup{:}

\medskip
\noindent
\textup{(i)}
\[
\left\{ \enspace
\begin{aligned}
  & \frac 1 p -\frac 1 q= \frac 4 {(m+2)(n+1)} \left( 1+ \frac m
  {2\mu} \right), \\
  & \frac 1 s= \frac{(m+2)(n-1)}4 \left( \frac 1 2 -\frac
  1 q\right) +\frac m {4 \mu}, \\
  & \gamma=\frac{n+1}{2} \left( \frac 1 2 -\frac 1 q \right)
  -\frac{m}{2\mu(m+2)},
\end{aligned}
\right.
\]
where $\mu \ge \mu_*$,
\[
\left\{ \enspace
\begin{aligned}
  & -\frac  1 {6\mu}<\gamma < \frac {47}{84} + \frac {25} {42 \mu}
  && \textup{if $n=2$, $m=1$,} \\
  & |\gamma- \gamma_\ast| < \gamma_d=
  \frac{2(2\mu-m)(n+1)}{\mu(m+2)(n-1)(2\mu_\ast-m)} && \textup{if
    $n\geq3$ or $n=2$, $m\geq2$,}
\end{aligned}
\right.
\]
and
\[
  \gamma_\ast=\frac 2 {m+2} +\frac m
        {2\mu(m+2)}-\frac{(2\mu-m)(n+1)}{2\mu(2\mu_\ast-m)}.
\]


\medskip
\noindent
\textup{(ii)} \ $n \ge 3$ or $n=2$, $m \ge 2$ and $r=2$,
\[
\left\{ \enspace
\begin{aligned}
  & \frac 1 s= \frac{(m+2)(n-1)}4 \left( \frac 1 2 -\frac 1 q\right)
  +\frac m {4 \mu}, \\
  &\gamma=\frac{n+1}{2} \left( \frac 1 2 -\frac 1 q\right) -\frac m
  {2\mu(m+2)},
\end{aligned}
\right.
\]
where $\mu \ge \max\{2,mn/2\}$ and
\[
  -\frac m {2\mu(m+2)} \le \gamma <\frac 3 {m+2} -\frac {n(2\mu
    -m)}{\mu (m+2)(n-1)}.
\]

\medskip
\noindent
\textup{(iii)} \ $n \ge 3$ or $n=2$, $m \ge 2$ and $s=2$,
\[
\left\{ \enspace
\begin{aligned}
  & \frac 1 r = 1-\frac m {4\mu} -\frac {(m+2)(n-1)} 4 \left( \frac 1
  p -\frac 1 2 \right), \\
  & \gamma=n \left(\frac 1 2 - \frac 1 q \right) - \frac 1 {m+2}, \\
\end{aligned}
\right.
\]
where $\mu \ge \max\{2,mn/2\}$ and
\[
  \frac {\mu(n+1)-mn} {\mu(m+2)(n-1)} <\gamma < \frac 2 {m+2} +\frac m
        {2\mu (m+2)}.
\]
\end{theorem}

\begin{remark}
We can rewrite the conditions of Theorem~\ref{4.5} in terms on $q$.

\smallskip
\noindent
(i) \ For $\mu \ge \mu_\ast$,
\begin{equation} \label{q1}
\begin{cases}
  \ds \frac8{63}\left(1-\frac4\mu\right) < \frac 1 q \le \frac 1 2 &
  \text{if $n=2, m=1$,} \\
  \ds\frac 1{p_2} < \frac 1 q +\frac 4 {(m+2)(n+1)} \left( 1 +\frac m
           {2\mu}\right)< \frac 1 {p_1} & \text{if $n\geq 3$ or $n=2$,
             $m\geq2$.}
\end{cases}
\end{equation}

\smallskip
\noindent
(ii) \ For $ \mu \ge \max\{2,mn/2\}$,
\begin{equation} \label{q2}
  \frac{2n}{\left(n+1\right)p_1} -\frac{n-1}{2(n+1)}
  -\frac{1}{(m+2)(n+1)} \left( 6+\frac m \mu \right) <\frac 1 q \le
  \frac 1 2.
\end{equation}

\smallskip
\noindent
(iii) \ For $ \mu \ge \max\{2,mn/2\}$,
\begin{equation} \label{q3}
  \frac 1 2 -\frac{1}{2(m+2)n} \left( 6+\frac m \mu \right) < \frac 1
  q < \frac 1 {q_1}.
\end{equation}
\end{remark}

\medskip

\begin{theorem}\label{thm4.5b}
Let $u$ solve the Cauchy problem \eqref{2.1} in the strip $S_T$. Then
\begin{equation} \label{4.32}
  \|u\|_{C_t^0\dot{H}_x^\gamma(S_T)} + \|u\|_{L^q(S_T)} \lesssim
  \|\varphi\|_{\dot{H}^\gamma(\R^n)}+ \|\psi\|_{\dot{H}^{\gamma-\frac
      2 {m+2}}(\R^n)} + \left\||D_x|^{\gamma - \gamma_0} f
  \right\|_{L^{p_0}(S_T)}
\end{equation}
provided that the exponents $p,\, q,\, r$, and $s$ satisfy \eqref{SI}
and $\mu \ge \max\{2,m/2\}$, $ q_0\le q<\infty$, where
\[
\gamma =\frac n 2 - \frac {n(m+2)+2}{q(m+2)}, \quad \gamma_0 =\frac 2
       {m+2} +\frac m {2\mu (m+2)} - \frac {n+1} 2\left( \frac 1 {p_0}
       -\frac 1 2 \right).
\]
\end{theorem}

\begin{corollary} \label{thm4.5c}
Under the conditions of Theorem \ref{thm4.5b}, one has
\begin{multline} \label{4.33}
  \|u\|_{C_t^0\dot{H}_x^\gamma(S_T)} + \|u\|_{L^q(S_T)} + \bigl\|
  |D_x|^{\gamma-\frac 1 {m+2}}u\bigr\|_{L^{q^\ast_0}(S_T)} \\ \lesssim
  \|\varphi\|_{\dot{H}^\gamma(\R^n)}+ \|\psi\|_{\dot{H}^{\gamma-\frac
      2 {m+2}}(\R^n)} + \bigl\| |D_x|^{\gamma- \frac 1 {m+2}} f
  \bigr\|_{L^{p_0^\ast}(S_T)},
\end{multline}
where $\gamma=\frac n 2 - \frac {(m+2)n+2}{q(m+2)}$ and $q_0^\ast \le
q <\infty$.
\end{corollary}
\begin{proof}
This follows by combining estimate \eqref{4.32} and
Remark~\ref{rem4.2} when $\mu=\mu_*$.
\end{proof}


An application of Theorem~\ref{thm4.5} yields:

\begin{corollary} \label{cor4.6}
Let $u$ solve the Cauchy problem
\begin{equation*}
\left\{ \enspace
\begin{aligned}
  &  \p_t^2 u -t^m \Delta u=f g  &&  \text{in $S_T$,} \\
  & u (0, \cdot)=\p_t u (0, \cdot) =0.
\end{aligned}
\right.
\end{equation*}
Then, for any $\mu \ge \mu_\ast$ and $0<R\leq\infty$,
\begin{equation}  \label{4.34}
  \|u\|_{C_t^0\dot{H}_x^\gamma(S_T\cap \Lambda_R)} + \|u\|_{L^s_t
    L^q_x (S_T \cap \Lambda_{R})}+ \|u\|_{L^\infty_t L^\delta_x(S_T
    \cap \Lambda_{R})}
  \lesssim \|f\|_{L^\sigma_t L^\rho_x(S_T \cap \Lambda_{R})}
  \|g\|_{L^s_t L^q_x(S_T \cap \Lambda_{R})},
\end{equation}
where $q$ is as in \eqref{q1},
\begin{align}
  \rho=\frac{\mu (m+2)(n+1)} {2(2\mu +m)}&, & \sigma &= \frac {\mu
    (n+1)}{ 2\mu -mn}, \label{4.38} \\
  \frac 1 s= \frac{(m+2)(n-1)}4 \left( \frac 1 2 -\frac 1 q\right)
  +\frac m {4 \mu}&, & \frac n \delta&=\frac n q + \frac 2 {m+2}
  \left( \frac 1 s -\frac m {4\mu} \right), \label{4.39}
\end{align}
and
\[
  \Lambda_{R}=\left\{(t,x) \in \R_+ \times \R^n \mid \, |x|+ \phi(t)<
  R \right\}.
\]
\end{corollary}


\begin{proof}
First we study the case $R=\infty$. Note that \eqref{4.39} gives that
\[
  n\left(\frac 1 2 -\frac 1 \delta\right)=\frac {n+1} 2\left(\frac 1 2
  -\frac 1 q\right) -\frac m {2\mu (m+2)}.
\]
Applying estimate \eqref{4.30} in case (i) together with the Sobolev
embedding $\dot{H}^{n\left(\frac 1 2 -\frac 1 \delta\right)}(\R^n)
\hookrightarrow L^\delta(\R^n )$, we have
\begin{equation*}
  \| u \|_{C_t^0\dot{H}_x^\gamma(S_T)} + \|u\|_{L^s_t L^q_x (S_T )} +
  \|u\|_{L^\infty_t L^\delta_x(S_T )} \lesssim \| f g\|_{L^r_t L^p_x
    (S_T)},
\end{equation*}
where $\frac 1 p=\frac 1 q +\frac 1 \rho$, $\frac 1 r=\frac 1 s +\frac
1 \sigma$. In addition, from H\"{o}lder's inequality,
\begin{equation}  \label{4.41}
 \| f g\|_{L^r_tL^p_x(S_T)}  \le  \|f\|_{L^\sigma_t
L^\rho_x(S_T)} \|g\|_{L^s_t L^q_x(S_T)}.
\end{equation}
Thus, estimate \eqref{4.34}  holds for $R=\infty$.

\smallskip

Now let $R<\infty$. Let $\chi$ denote the characteristic function of
$S_T \cap \Lambda_{R}$. If $u$ solves $\p_t^2 u -t^m \Delta u= f g$
with vanishing initial data and $u_\chi$ solves $\p_t^2 u_\chi -t^m
\Delta u_\chi = \chi f g$ with vanishing initial data, then $u=u_\chi$
in $S_T \cap \Lambda_{R}$ due to finite propagation speed (see
\cite{Taniguchi}). Therefore,
\begin{multline*}
  \| u \|_{C_t^0\dot{H}_x^\gamma(S_T\cap \Lambda_R)} + \|u\|_{L^s_t L^q_x (S_T \cap
    \Lambda_{R})} + \| u \|_{L^\infty_t L^\delta_x(S_T \cap
    \Lambda_{R})} \\ = \|u_\chi\|_{C_t^0\dot{H}_x^\gamma(S_T)} + \|
  u_\chi\|_{L^s_t L^q_x (S_T )} + \|u_\chi\|_{L^\infty_t
    L^\delta_x(S_T )} \le \|\chi f\|_{L^\sigma_t L^\rho_x(S_T)} \|\chi
  g\|_{L^s_t L^q_x(S_T)}.
\end{multline*}
Consequently, estimate \eqref{4.34} holds.
\end{proof}

As another application of Theorem~\ref{thm4.5} we have:
\begin{corollary} \label{cor4.10}
Let $u$ be a solution of
\begin{equation*}
\left\{ \enspace
\begin{aligned}
  & \p_t^2 u-t^m \Delta u=F(v) && \text{in } S_T, \\
  & u(0, \cdot)=\p_t u (0, \cdot) =0.
\end{aligned}
\right.
\end{equation*}
If $q<\infty$ and $\frac 1 {m+2} \le \gamma =\frac n 2 -
\frac {n(m+2)+2}{q(m+2)}\le \frac{m+3}{m+2}$, then
\begin{multline}  \label{4.40}
 \|u\|_{C_t^0\dot{H}_x^\gamma(S_T)} + \|u\|_{ L^q (S_T )}  + \|
 |D_x|^{\gamma-\frac 1 {m+2}} u\|_{ L^{q_0^\ast}(S_T )} \\ \lesssim
 \|F'(v)\|_{L^{ \frac{\mu_\ast} 2 }(S_T )} \bigl\||D_x|^{\gamma-\frac 1
   {m+2}} v\bigr\|_{ L^{q_0^\ast}(S_T)}.
\end{multline}
\end{corollary}

\begin{proof}
This follows from estimate \eqref{4.33} by taking fractional
derivatives. Indeed, for $0 \le \gamma -\frac 1 {m+2} \le~1$, one has
\begin{multline*}
  \|u\|_{C_t^0\dot{H}_x^\gamma(S_T)} + \|u\|_{ L^q (S_T )} + \bigl\|
  |D_x|^{\gamma-\frac 1 {m+2}} u \bigr\|_{ L^{q_0^\ast}(S_T )}
  \\ \lesssim \bigl\| |D_x|^{\gamma-\frac 1 {m+2} } \left( F(v)
  \right)\bigr\|_{ L^{ p_0^\ast }(S_T) } \lesssim \|F'(v)\|_{L^{
      \frac{\mu_\ast} 2 }(S_T )} \bigl\||D_x|^{\gamma-\frac 1 {m+2}} v\bigr\|_{
    L^{q_0^\ast}(S_T)}. 
\end{multline*}
\end{proof}


\section{Solvability of the semilinear generalized Tricomi equation}\label{sec5}

In this section, we will apply Theorems~\ref{thm4.5} and \ref{thm4.5b}
and Corollaries~\ref{thm4.5c} to \ref{cor4.10} with $\mu=\mu_\ast$ to
establish the existence and uniqueness of the solution $u$ of problem
\eqref{1.1}. Thereby, we will use the following iteration scheme: For
$j \in \mathbb N_0$, let $u_j$ be the solution of
\begin{equation} \label{5.1}
\left\{ \enspace
\begin{aligned}
 & \p_t^2 u_{j} -t^m \Delta u_{j} = F(u_{j-1}) && \text{in } \R_+
  \times \R^n, \\
 & u_{j}(0, \cdot)=\varphi, \enspace \p_t u_{j}(0,\cdot)=\psi,
\end{aligned}
\right.
\end{equation}
where $u_{-1}= 0$.
Notice that, for $\mu=\mu_\ast$, the exponents from \eqref{4.30}
in case~(i) are
\[
\gamma_\ast =\frac1{m+2}, \quad \gamma_d = \frac
      {2\left(n+1\right)}{\mu_\ast (m+2)(n-1)}.
\]
In order to get the existence of solutions of the Cauchy problem
\eqref{1.1} as stated in Theorems \ref{thm1.1}, \ref{thm1.5}, and
\ref{thm1.6}, we need to show that, for the sequences
$\{u_j\}_{j=0}^\infty$ and $\{F(u_j)\}_{j=0}^\infty$ defined by
\eqref{5.1}, there exist a $T>0$ and a function $u$ such that
\begin{align}
  u_j &\rightarrow u &&\hspace{-40mm}\text{in } L_{\text{loc}}^1(S_T) \quad \text{as
  } j \rightarrow \infty, \label{5.2} \\
  F(u_j) &\rightarrow F(u) &&\hspace{-40mm}\text{in } L_{\text{loc}}^1(S_T) \quad \text{as
  } j \rightarrow \infty.  \label{5.3}
\end{align}
From \eqref{5.2} and \eqref{5.3}, one obviously has that the limit
function $u$ solves problem \eqref{1.1} in $S_T$.

Furthermore, let $u,\, \tilde u$ both solve the Cauchy problem
\eqref{1.1} in $S_T$. Then $v=u-\tilde u$ satisfies
\begin{equation} \label{U1}
\left\{ \enspace
\begin{aligned}
& \p_t^2 v -t^m \Delta v  = G(u, \tilde u) v  && \text{in }  S_T, \\
& v (0, \cdot)=  \p_t  v (0, \cdot)=0,
\end{aligned}
\right.
\end{equation}
where $G(u, \tilde u) =\frac{F(u) -F(\tilde u)}{u -\tilde u}$ if
$u\neq \tilde u$ and $G(u,u)=F'(u)$. For certain $s,\, q \ge 2$, we
will show that $v \in L^s_t L^q_x(S_T)$ and then
\begin{equation} \label{U2}
 \|v\|_{L^s_t L^q_x(S_T)} \le \frac 1 2\, \|v\|_{L^s_t L^q_x(S_T)}.
\end{equation}
Uniqueness of the solution of the Cauchy problem \eqref{1.1} in $S_T$
follows.



\subsection{Proof of Theorem~\ref{thm1.1}}


\subsubsection{Case $\kappa_1 < \kappa < \kappa_\ast$}\label{case_1}

From the assumptions of Theorem~\ref{thm1.1}, we have
\[
\gamma=\frac {n+1} 4 -\frac{n+1}{\mu_\ast (\kappa-1)}-\frac m {2
  \mu_\ast (m+2)}
\]
and
\begin{equation}   \label{5.4}
  q=\frac {\mu_\ast\left(\kappa-1\right)}2 , \quad \frac 1 s=
  \frac{(m+2)(n-1)}4 \left( \frac 1 2 -\frac 1 q\right) +\frac m {4
    \mu_\ast}.
 \end{equation}
Thus,
\[
  \gamma=\frac {n+1} 2 \left( \frac 1 2 -\frac 1 q \right) -\frac m {
    2 \mu_\ast (m+2)}, \quad \frac 1 {m+2} - \frac {2(n+1)}{ \mu_\ast
    (m+2)(n-1)} < \gamma < \frac 1 {m+2}.
\]

\paragraph{\bf Existence.} \ In order to show \eqref{5.2}, set
\begin{equation}  \label{5.5}
   H_j(T)= \|u_j\|_{C_t^0\dot{H}_x^\gamma(S_T)} + \|u_j\|_{L^s_t
     L^q_x(S_T)}, \quad N_j(T)=\|u_j - u_{j-1}\|_{L^s_t L^q_x(S_T)}.
\end{equation}

We claim that there exists a constant $\varepsilon_0 >0$ small
such that
\begin{equation} \label{5.7}
2 T^{\frac 1 q -\frac 1 s} H_0(T)  \le \varepsilon_0
\end{equation}
and
\begin{equation} \label{5.6}
H_j (T) \le 2  H_0(T), \quad N_j(T) \le \frac { 1} { 2}\, N_{j-1}(T).
\end{equation}

Indeed, from the iteration scheme \eqref{5.1}, we have
\begin{equation}  \label{5.8}
  \left( \p_t^2 - t^m \Delta \right) (u_{j+1} -u_{k+1}) = G(u_j, u_k)
  (u_j -u_k).
\end{equation}
Note that in \eqref{4.38}
\[
  \rho=\sigma=\frac {\mu_\ast} 2
\]
when $\mu=\mu_\ast$. Thus, from \eqref{4.34} and condition
\eqref{1.2},
\begin{multline}  \label{5.9}
  \|u_{j+1} -u_{k+1}\|_{C_t^0\dot{H}_x^\gamma(S_T)} + \|u_{j+1} -u_{k+1}\|_{L^s_t L^q_x(S_T)}\\
\begin{aligned}
  & \lesssim \|G(u_j, u_k)\|_{L^{\frac {\mu_\ast} 2}(S_T)} \|u_{j}
  -u_{k}\|_{L^s_t L^q_x(S_T)} \\
  & \lesssim\bigl( \|u_j\|_{L^{q} (S_T)}^{\kappa-1} + \|u_k\|_{L^{q}
    (S_T)}^{\kappa-1} \bigr) \|u_{j} -u_{k}\|_{L^s_t L^q_x(S_T)}.
\end{aligned}
\end{multline}
Note that $s > q$ for $\kappa < \kappa_\ast$.  By H\"{o}lder's
inequality, we arrive at
\begin{equation} \label{5.10}
  \|u_j\|_{L^{q} (S_T)}\le T^{\frac 1 q -\frac 1 s} \|u_j\|_{L^s_t
    L^q_x(S_T)}.
\end{equation}
Since $u_{-1}= 0$, \eqref{5.9} together with \eqref{5.10} implies
that
\[
  \| u_{j+1} -u_0\|_{L^s_t L^q_x(S_T)} +
  \|u_{j+1}-u_0\|_{C_t^0\dot{H}_x^\gamma(S_T)} \lesssim
  T^{(\kappa-1)\left(\frac 1 q -\frac 1 s\right)}
  \|u_j\|^\kappa_{L^s_t L^q_x (S_T)}.
\]
From the Minkowski inequality, we have that there exists an $\ve_0$
with $0<\ve_0\leq 2^{-2/\left(\kappa-1\right)}$ such that
\[
  H_{j+1}(T) \le H_0(T) + \frac 1 2\, H_j(T) \quad \text{if } T^{\frac
    1 q -\frac 1 s} H_j(T) \le \varepsilon_0.
\]
Therefore, by induction on $j$,
\begin{equation}   \label{5.11}
  H_j(T) \le 2 H_0(T) \quad \text{if } 2T^{\frac 1 q -\frac 1 s}
  H_0(T) \le \varepsilon_0.
\end{equation}
Taking $k=j-1$ in \eqref{5.8}, estimates \eqref{5.9} to \eqref{5.11}
yield that
\[
  N_{j+1}(T) \le \frac 1 2\, N_j(T) \quad \text{if } 2 H_0(T) T^{\frac
    1 q -\frac 1 s} \le \varepsilon_0,
\]
which together with \eqref{5.11} implies that \eqref{5.6} holds as
long as \eqref{5.7} holds.

Since $u_{-1}\equiv 0$ and $u_0$ is a solution of problem \eqref{2.2},
we have from \eqref{4.14} that, for $\varphi \in \dot{H}^\gamma(\R^n)$
and $\psi \in \dot{H}^{\gamma - \frac 2 {m+2}}(\R^n)$,
\[
  N_0(T)\le H_0(T)\lesssim \|\varphi\|_{\dot{H}^\gamma(\R^n)} +
  \|\psi\|_{\dot{H}^{\gamma - \frac 2 {m+2}}(\R^n)}.
\]
Thus, by choosing $T>0$ small, \eqref{5.7} holds.  Consequently, there
is a function $u \in C_t^0\dot{H}_x^\gamma(S_T)\cap \,
L^s_tL^q_x(S_T)$ such that
\begin{equation}   \label{5.12}
  u_j \rightarrow u \quad \text{in } L^s_t L^q_x(S_T) \quad \text{as }
  j \rightarrow \infty,
\end{equation}
and, therefore, \eqref{5.2} holds. It also follows that $u_j$
converges to $u$ almost where. By Fatou's lemma, it follows that
\begin{equation} \label{H1}
  \|u\|_{C_t^0\dot{H}_x^\gamma(S_T)} + \|u\|_{L^s_t L^q_x(S_T)}
  \le \liminf \limits_{ j \rightarrow \infty}
  \left(\|u_j\|_{C_t^0\dot{H}_x^\gamma(S_T)} + \|u_j\|_{L^s_t
    L^q_x(S_T)} \right) \le 2 H_0(T),
\end{equation}
which shows that estimate \eqref{1.4} holds.

\smallskip

Now we prove \eqref{5.3}. It suffices to show that $F(u)$ is bounded
in $L^r_t L^p_x(S_T)$ and $F(u_j)$ converges to $F(u)$ in $L^r_t
L^p_x(S_T)$ as $j \rightarrow \infty$, where $p=q/\kappa$ and
$\frac1r= 1- \frac{m}{4\mu_*} - \frac{(m+2)(n-1)}4
\bigl(\frac1p-\frac12\bigr)$. In fact, $r \kappa < s$ if $ \kappa
< \kappa_\ast$, thus, for $q= p \kappa$, by condition \eqref{1.2}
and H\"older's inequality, we have
\[
  \|F(u)\|_{L^r_t L^p_x (S_T)}\lesssim \|u\|_{L^{r \kappa}_t L^{p
      \kappa}_x (S_T)}^\kappa\lesssim T^{\frac 1 {r } -\frac \kappa s}
  \|u\|_{L^s_t L^q_x(S_T)}^\kappa.
\]
Moreover, in view of $\frac 1 p -\frac 1 q =\frac 1 r -\frac 1 s = \frac
2{\mu_\ast},$ by H\"older's inequality and estimates
\eqref{5.9}-\eqref{5.11} and \eqref{H1}, we have
\begin{multline*}
  \| F(u_j) -F(u)\|_{L^r_t L^p_x (S_T)} \le \|G(u_j,
  u)\|_{L^{\mu_\ast/2}(S_T)} \|u_j -u\|_{L^s_t L^q_x (S_T)}\\ \lesssim
  T^{(\kappa-1) \bigl( \frac 1 q -\frac 1 s \bigr)} H_0(T)^{\kappa-1}
  \|u_j -u\|_{L^s_t L^q_x (S_T)} \lesssim \|u_j -u\|_{L^s_t L^q_x
    (S_T)}.
\end{multline*}
Applying \eqref{5.12}, we have that $F(u_j)$ converges to $F(u)$ in
$L^r_t L^p_x(S_T)$ and, therefore, \eqref{5.3} holds.

\smallskip

From \eqref{5.2} and \eqref{5.3}, we have that the limit function $u
\in C^0_t \dot{ H}^\gamma_x(S_T)) \cap L^s_tL^q_x(S_T) $ solves the
Cauchy problem \eqref{1.1} in $S_T$.

\paragraph{\bf Uniqueness} \
Suppose $u,\,\tilde u \in C([0, T],\dot{ H}^\gamma(\R^n)) \cap
L^s_tL^q_x(S_T)$ solve the Cauchy problem \eqref{1.1} in $S_T$. Then
$v=u - \tilde u \in C([0, T],\dot{ H}^\gamma(\R^n)) \cap
L^s_tL^q_x(S_T)$ is a solution of problem \eqref{U1}. From
Corollary~\ref{cor4.6}, we have that
\begin{align*}
  \| v \|_{L^s_t L^q_x(S_T)} & \le C  \bigl( \|u\|_{L^{q}
    (S_T)}^{\kappa-1} + \|\tilde u\|_{L^{q} (S_T)}^{\kappa-1} \bigr)
  \| v\|_{L^s_t L^q_x(S_T)} && \textup{(by \eqref{4.34} and
    \eqref{1.2})}\\
  & \le C  T^{(\kappa-1)(\frac 1 q -\frac 1 s)} \bigl( \|u\|_{L^s_t
    L^q_x (S_T)}^{\kappa-1} + \|\tilde u\|_{L^s_t L^q_x
    (S_T)}^{\kappa-1} \bigr) \| v \|_{L^s_t L^q_x(S_T)} && \textup{(by
    H\"{o}lder's inequality)} \\
  & \le C 2^\kappa \bigl( T^{\frac 1 q -\frac 1 s}
  H_0(T)\bigr)^{\kappa-1} \|v \|_{L^s_t L^q_x(S_T)} && \textup{(by
    \eqref{H1})} \\
  & \le  \frac 1 2 \| v\|_{L^s_t L^q_x(S_T)}. && \textup{(by
    \eqref{5.7})}
\end{align*}
Thus  \eqref{U2} holds and $u=\tilde u$ in $S_T$.


\subsubsection{Case $\kappa_\ast\leq \kappa$ if $n=2$
or $\kappa_\ast \le \kappa \leq \kappa_3$ if $n \ge 3$.}\label{case_2}

${ }$
\smallskip

\paragraph{\bf Existence} \
From the assumptions of Theorem~\ref{thm1.1}, we have
\[
  \gamma=\frac {n} 2 -\frac{4}{ (m+2) (\kappa-1)}, \quad
  s=q=\frac{\mu_\ast\left(\kappa-1\right)} 2.
\]
Thus,
\[
  \frac 1 {m+2} \le \gamma = \frac n 2 - \frac{(m+2)n+2}{q(m+2)} \le
  \frac{m+3}{m+2}.
\]

To show  \eqref{5.2}, we set
\begin{align*}
  H_j(T)= \|u_j\|_{C_t^0\dot{H}_x^\gamma(S_T)} + \|u_j\|_{L^q(S_T)} +
  \| |D_x|^{\gamma -\frac 1 {m+2}} u_j\|_{L^{q_0^\ast}(S_T)},
\end{align*}
and
\begin{equation} \label{5.14}
  N_j(T)=\|u_j - u_{j-1}\|_{L^{q_0^\ast}(S_T \cap \Lambda_{R})}.
\end{equation}

We claim that there exists a constant $\varepsilon_0>0$ such that
\begin{equation} \label{5.16}
  H_0(T)  \le \varepsilon_0,
\end{equation}
and
\begin{equation} \label{5.15}
  H_j(T) \le 2  H_0(T), \quad N_j(T) \le \frac { 1} { 2} N_{j-1}(T).
\end{equation}
Indeed, since $u_{-1}= 0$, from the iteration scheme \eqref{5.1}, we
have
\begin{equation}  \label{5.17}
  \bigl( \p_t^2 - t^m \Delta \bigr) (u_{j+1} -u_{0}) = F(u_j).
\end{equation}
Thus, estimate \eqref{4.40} together with condition \eqref{1.2} yields
that, for $0 \le \gamma -\frac 1 {m+2} \le 1$,
\begin{align*}
  H_{j+1} (T) & \le H_0(T) + C \|F'(u_j)\|_{L^{\frac {\mu_\ast}
      2}(S_T)} \bigl\| |D_x|^{\gamma -\frac 1 {m+2}} u_{j}
  \bigr\|_{L^{q_0^\ast}(S_T)} \\
  & \le H_0(T) + C \|u_j\|_{L^{q} (S_T)}^{\kappa-1} \bigl\| |D_x|^{\gamma
    -\frac 1 {m+2}} u_{j} \bigr\|_{L^{q_0^\ast}(S_T)} \\
  & \le H_0(T) +C H_j(T)^\kappa.
\end{align*}
Therefore, by induction, we have that
\[
  H_j(T) \le 2 H_0(T) \quad
  \text{if } C 2^\kappa H_0(T)^{\kappa-1}<1.
\]
Consequently,
\begin{equation}\label{5.18}
  H_j(T) \le 2 H_0(T) \quad \text{if } H_0(T)\le \varepsilon_0
\end{equation}
for some $\varepsilon_0 >0$ small. Notice that, for $q$ and $s$ from
\eqref{5.4}, when $q=s$, so $q=s=q_0^\ast.$ Hence, by using
estimates \eqref{5.9}-\eqref{5.11} together with
\eqref{5.18}, we get that for $N_j$ defined in \eqref{5.14},
\begin{equation} \label{5.19}
  N_j(T) \le \frac 1 2 N_{j-1}(T) \quad \text{if } H_0(T)\le
  \varepsilon_0.
\end{equation}
Estimates \eqref{5.18} and \eqref{5.19} tell us that \eqref{5.15}
holds as long as \eqref{5.16} holds. To get \eqref{5.16}, from
estimate \eqref{4.33} (with $f=0$) we have that, for $\varphi \in
\dot{H}^\gamma(\R^n)$ and $\psi \in \dot{H}^{\gamma - \frac 2
  {m+2}}(\R^n)$,
\begin{equation}  \label{5.20}
  H_0(T) \lesssim \|\varphi\|_{\dot{H}^\gamma(\R^n)} +
  \|\psi\|_{\dot{H}^{\gamma - \frac 2 {m+2}}(\R^n)}.
\end{equation}
Due to the continuity of the norm in $L^q(S_T)$, \eqref{5.16} holds
for some $T>0$ small. (If $ \|\varphi\|_{\dot{H}^\gamma(\R^n)} +
\|\psi\|_{\dot{H}^{\gamma - \frac 2 {m+2}}(\R^n)}$ is small, then
\eqref{5.16} holds for any $T>0$, consequently, we get global
existence.)

Note that $q = \mu_\ast(\kappa-1)/2 \ge q_0^\ast$ when $\kappa \ge
\kappa_\ast$. Thus, from H\"older's inequality and \eqref{5.20},
\begin{equation} \label{5.21}
  N_0(T)= \|u_0\|_{L^{q_0^\ast}(S_{T} \cap \Lambda_{R})} \lesssim
  \|u_0\|_{L^{q}(S_{T})} \lesssim H_0(T).
\end{equation}
From estimates \eqref{5.16}, \eqref{5.15}, and \eqref{5.21}, we get
that there exists a function $u \in C_t^0\dot{ H}_x^\gamma(S_T)
\cap $ $ L^q(S_{T})$ with $|D_x|^{\gamma-\frac 1 {m+2}} u \in
L^{q_0^\ast}(S_T)$ such that
\begin{equation} \label{5.22}
  u_j \rightarrow u \quad \text{in } \, L^{q_0^\ast}(S_T \cap
  \Lambda_{R}) \quad \text{as } \, j \rightarrow \infty,
\end{equation}
and \eqref{5.2} holds. Thus, from Fatou's lemma and \eqref{5.15}, it
follows that
\begin{equation} \label{H2}
  \|u\|_{C_t^0\dot{H}^\gamma_x(S_T)} + \|u\|_{L^q(S_T)} + \bigl\|
  |D_x|^{\gamma -\frac 1 {m+2}} u\bigr\|_{L^{q_0^\ast}(S_T)} \le 2
  H_0(T)
\end{equation}
and  $u$ satisfies estimate \eqref{1.4}.

Since $q=\mu_\ast(\kappa-1)/2 \geq \kappa$ when $\kappa \ge \kappa_\ast$,
we have from condition \eqref{1.2} that $F(u)$ is locally integrable
for $u \in L^q(S_T)$. By H\"older's inequality,
\begin{multline*}
  \int_{S_T \cap \Lambda_{R}} | F(u_j) -F(u)|\, dt dx = \int_{S_T \cap
    \Lambda_{R}} |G(u_j, u)| \left| u_j -u\right| dt dx \\
  \le \| G(u_j,u)\|_{L^{p_0^\ast}(S_T \cap \Lambda_{R})} \|u_j
  -u\|_{L^{q_0^\ast}(S_T \cap \Lambda_{R})}.
\end{multline*}
Note that $ p_0^\ast < \mu_\ast/2$. Thus, from condition  \eqref{1.2} we have that
\begin{multline*}
  \| G(u_j, u)\|_{L^{p_0^\ast}(S_T \cap \Lambda_{R})} \lesssim
  \|u_j\|_{L^{p_0^\ast (\kappa-1)} (S_T \cap \Lambda_{R})}^{\kappa -1}
  + \|u\|_{L^{p_0^\ast (\kappa-1)}(S_T \cap \Lambda_{R})}^{\kappa -1}
  \\
  \lesssim \|u_j\|^{\kappa -1}_{L^q(S_T \cap \Lambda_{R})} +
  \|u\|^{\kappa -1}_{L^q(S_T \cap \Lambda_{R})} \lesssim
  H_0({T})^{\kappa-1},
\end{multline*}
which together with \eqref{5.22} implies that $F(u_j) \rightarrow F(u)$ in
$L^1_{\text{loc}}(S_T)$. Hence, \eqref{5.3} holds.

From \eqref{5.2} and \eqref{5.3}, we have that the limit function $u
\in C_t^0\dot{H}_x^\gamma(S_T) \cap L^q(S_T )$ with $|D_x|^{\gamma-\frac
  1 {m+2}} u \in L^{q_0^\ast}(S_T )$ is a weak solution of Cauchy
problem \eqref{1.1} in $S_T$.

\smallskip

\paragraph{\bf Uniqueness}  \
Suppose $u, \tilde u \in C_t^0\dot{H}_x^\gamma(S_T) \cap L^q (S_T)$
with $|D_x|^{\gamma-\frac 1 {m+2}} u,\, |D_x|^{\gamma-\frac 1 {m+2}}
\tilde u \in L^{q_0^\ast}(S_T )$ solve the Cauchy problem \eqref{1.1}
in $S_T$.  Then $v=u-\tilde u \in C_t^0\dot{H}_x^\gamma(S_T) \cap
L^q(S_T) $ is a weak solution of problem \eqref{U1}. Thus, it follows
from Corollary~\ref{cor4.6} that
\[
\begin{aligned}
  \| v \|_{L^q (S_T)} & \le C \left( \|u\|_{L^{q} (S_T)}^{\kappa-1} +
  \|\tilde u \|_{L^{q} (S_T)}^{\kappa-1} \right) \| v \|_{ L^q(S_T)}
  && \textup{(by \eqref{4.34} and \eqref{1.2})}\\
  & \le C 2^\kappa H_0(T)^{\kappa-1} \|v\|_{ L^q(S_T)} && \textup{(by
    \eqref{H2})}\\
  & \le \frac 1 2\, \| v \|_{L^s_t L^q_x(S_T)} && \textup{(by
    \eqref{5.16}).}
\end{aligned}
\]
Thus \eqref{U2} holds and  $u=\tilde u$ in $S_T$.


\subsubsection{Case $n\geq3$ and $\kappa > \kappa_3$, $\kappa \in
   \mathbb N$}\label{case_3}

${ }$
\smallskip

\paragraph{\bf Existence}  \
From the assumptions of Theorem~\ref{thm1.1}, we have
\[
  \gamma=\frac {n} 2-\frac{4}{ (m+2) (\kappa-1)}, \quad
  s=q=\frac{\mu_\ast\left(\kappa-1\right) } 2, \quad F(u)=\pm u^\kappa,
\]
and
\[
  \gamma= \frac n 2 - \frac {(m+2)n+2}{q(m+2)} > 1 + \frac{1}{m+2}.
\]

To verify \eqref{5.2}, we set
\[
  H_j(T)= \|u_j\|_{C_t^0\dot{H}_x^\gamma(S_T)} + \sup_{q_0^\ast \le
    \tau \le \frac{\mu_\ast\left(\kappa-1\right)} 2} \bigl\| |D_x|^{\frac
    {(m+2)n+2} {\tau(m+2)} -\frac 4 {(m+2)(\kappa-1)} }
  u_j\bigr\|_{L^\tau(S_T)}
\]
and
\[
  N_j(T)=\|u_j - u_{j-1}\|_{L^{q_0^\ast}(S_T \cap \Lambda_{R})}.
\]

We claim that there exists a constant $\varepsilon_0>0$ such that
\begin{equation}  \label{5.24}
 H_0(T)  \le \varepsilon_0
 \end{equation}
 and
\begin{equation}  \label{5.23}
H_j(T) \le 2  H_0(T), \quad N_j(T) \le \frac { 1} { 2} N_{j-1}(T).
\end{equation}

In fact, applying Minkowski's inequality and estimate \eqref{4.33}  (with
$\varphi = \psi  = 0$),
\begin{equation} \label{5.25}
  H_{j+1} (T)\le H_0(T) + C \sup \limits_{q_0^\ast \le \tau \le
   \mu_\ast (\kappa-1)/2} \| |D_x|^{ \frac n 2 -\frac 1
    {m+2}-\frac 4 {(m+2)(\kappa-1)}} \bigl( u_{j}^\kappa \bigr)
  \|_{L^{p_0^\ast}(S_T)}.
\end{equation}
Note that $\alpha =\frac n 2 -\frac 1 {m+2}-\frac 4
{(m+2)(\kappa-1)}>1$ when $\kappa > \kappa_3$. Thus, $|D_x|^\alpha
\bigl( u_j^\kappa \bigr)$ can be expressed as a finite linear
combination of $\prod \limits_{\ell =1}^\kappa |D_x|^{\alpha_\ell}
u_j$, where $ 0\le \alpha_\ell \le \alpha$ ($1\le \ell \le \kappa$)
and $ \sum \limits_{\ell =1}^\kappa \alpha_\ell = \alpha. $ By
H\"older's inequality, $\bigl\| |D_x|^{ \alpha} \bigl( u_{j}^\kappa
\bigr) \bigr\|_{L^{p_0^\ast}(S_T)} $ is dominated by a finite sum of
terms of the form $\prod \limits_{\ell =1}^\kappa
\bigl\||D_x|^{\alpha_\ell} u_j\bigr\|_{L^{\tau_\ell}(S_T)}$, where $\sum
\limits_{\ell =1}^\kappa 1/{\tau_\ell} = 1/{p_0^\ast}$. We
choose $\tau_\ell$ so that
\[
  \alpha_\ell=\frac {n(m+2)+2}{\tau_\ell (m+2)}-\frac 4
        {(m+2)(\kappa-1)}.
\]
Then
\[
  q_0^\ast \le \tau_\ell \le \frac {\mu_\ast\left(\kappa-1\right)}2,
  \quad \sum \limits_{\ell =1}^\kappa \frac 1 {\tau_\ell} =\frac 1
        {p_0^\ast},
\]
and, therefore,
\[
  \bigl\||D_x|^{\alpha_\ell} u_j\bigr\|_{L^{\tau_\ell}(S_T)} \le H_j(T),
\]
which together with \eqref{5.25} yields that
\[
  H_{j+1}(T) \le  H_0(T)+ C_\kappa H_j(T)^\kappa.
\]
By induction, we have that
\begin{equation} \label{5.26}
H_j(T) \le 2 H_0(T) \quad \text{if }   H_0(T)\le
\varepsilon_0.
\end{equation}
For $q$ and $s$ from \eqref{5.4}, when $q=s$, so $q=s=q_0^\ast$.
Hence, by estimates \eqref{5.9}-\eqref{5.11} and together with
\eqref{5.26}, we get that
\begin{equation} \label{5.27}
  N_j(T) \le \frac 1 2 N_{j-1}(T) \quad \text{if } H_0(T)\le
  \varepsilon_0.
\end{equation}
From \eqref{5.26} and \eqref{5.27}, we get that \eqref{5.23} holds as
long as \eqref{5.24} holds.

Note that
\begin{equation} \label{5.28}
  \frac {n(m+2)+2}{\tau(m+2)} -\frac 4 {(m+2)(\kappa -1)}=0,
\end{equation}
for $\tau= \mu_\ast(\kappa -1)/2$ and
\begin{equation} \label{5.29}
  \frac {n(m+2)+2}{\tau(m+2)} -\frac 4 {(m+2)(\kappa -1)}= \gamma
  -\frac 1 {m+2}.
\end{equation}
for $\tau=q_0^\ast$. On the other hand, we have from \eqref{4.33}
(with $f=0$) that, for $\varphi \in \dot{H}^\gamma(\R^n)$ and $\psi
\in \dot{H}^{\gamma - \frac 2 {m+2}}(\R^n)$,
\begin{multline}  \label{5.30}
   \|u_0\|_{C_t^0\dot{H}^\gamma_x(S_T)} +\|u_0\|_{L^{\frac{ \mu_\ast
         (\kappa -1)} 2}(S_T)} + \bigl\| |D_x|^{\gamma -\frac{1}{m+2}}
   u_0 \bigr\|_{L^{p_0^\ast}(S_T)} \\
   \lesssim \|\varphi\|_{\dot{H}^\gamma(\R^n)} +
   \|\psi\|_{\dot{H}^{\gamma - \frac 2 {m+2}}(\R^n)}.
\end{multline}
By interpolation together with \eqref{5.28}-\eqref{5.30}, we conclude
that
\[
  H_0(T) \lesssim \|\varphi\|_{\dot{H}^\gamma(\R^n)} +
  \|\psi\|_{\dot{H}^{\gamma - \frac 2 {m+2}}(\R^n)}.
\]
It follows that \eqref{5.24} holds by choosing $T >0$ small. (We can
take $T =\infty$ if $\|\varphi\|_{\dot{H}^\gamma(\R^n)} +
\|\psi\|_{\dot{H}^{\gamma - \frac 2 {m+2}}(\R^n)}$ is small which then
yields global existence.)

From H\"older's inequality and \eqref{5.28},
\begin{equation}  \label{5.31}
  N_0(T)= \|u_0\|_{L^{q_0^\ast}(S_T \cap \Lambda_{R})} \le C_R
  \left\|u_0\right\|_{L^{\frac{\mu_\ast (\kappa -1)} 2}(S_T)} \le C_R H_0(T)
  <\infty.
\end{equation}
Therefore, we have from \eqref{5.23}, \eqref{5.24}, and \eqref{5.31}
that there exists a function $u \in C_t^0\dot{H}_x^\gamma(S_T) \cap
L^q(S_T)$ with $|D_x|^{\gamma-\frac 1 {m+2}} u \in L^{q_0^\ast}(S_T)$
such that
\[
  u_j \rightarrow u \quad \text{in} \ \ L^{q_0^\ast}(S_T \cap
  \Lambda_{R}) \quad \text{as } j \rightarrow \infty,
\]
and, therefore, \eqref{5.2} holds. Thus, from Fatou's lemma and
\eqref{5.23},
\begin{equation} \label{H3}
   \|u\|_{C_t^0\dot{H}_x^\gamma(S_T)} + \|u\|_{L^q(S_T)} + \bigl\|
   |D_x|^{\gamma -\frac 1 {m+2}} u \bigr\|_{L^{q_0^\ast}(S_T)} \le 2
   H_0(T)
\end{equation}
and $u$ satisfies estimate \eqref{1.4}.

Note that $q=\mu_\ast(\kappa-1)/2 \ge \kappa$ when $\kappa >
\kappa_3$. Thus, for $u \in L^q(S_T)$, by H\"older's inequality and
condition \eqref{1.2}, we get that $F(u)$ is locally integrable and
$F(u_j)$ convergences to $F(u)$ in $L^1_{\textup{loc}}(S_T)$, and
hence \eqref{5.3} holds.

Applying \eqref{5.2}, \eqref{5.3}, it follows that the limit
function $u \in C_t^0\dot{H}_x^\gamma(S_T)\cap L^q (S_T)$ with
$|D_x|^{\gamma-\frac 1 {m+2}} u \in L^{q_0^\ast}(S_T) $ is a weak
solution of the Cauchy problem \eqref{1.1} in $S_T$.
\smallskip

\paragraph{\bf Uniqueness} \
This follows from the same arguments as in \ref{case_2}.
\qed


\subsection{Proof of Theorem~\ref{thm1.5}}

From the assumption of Theorem~\ref{thm1.5}, we have
\begin{align*}
  \gamma &= \frac {n} 2 -\frac{4}{(m+2) (\kappa-1)}, \\
  \frac1q &= \frac 1 {(m+2)(n+1)} \left(\frac 8 {\kappa-1} -\frac m
          {\mu_\ast} \right) -\frac{n-1}{2(n+1)}, \\
\intertext{and}
  \frac 1 s&= \frac{(m+2)(n-1)}4 \left( \frac 1 2 -\frac 1 q\right)
  +\frac m {4 \mu_\ast}.
\end{align*}
Thus,
\[
  \gamma=\left( \frac {n+1} 2\right) \left( \frac 1 2 -\frac 1 q
  \right) -\frac m { 2 \mu_\ast (m+2)}
\]
and
\[
  \frac 1 {m+2} \le \gamma <
  \frac 1 {m+2} + \frac {2(n+1)}{\mu_\ast (m+2)(n-1)},
\]
where $\kappa_\ast \le \kappa < \kappa_2$.

\smallskip

To show \eqref{5.2}, we set
\[
  H_j(T)= \|u_j\|_{C_t^0\dot{H}_x^\gamma(S_T)} + \|u_j\|_{L^s_t
    L^q_x(S_T)} + \|u_j -u_0\|_{L^\infty_t L^\delta_x(S_T)}
\]
and
\[
  N_j(T)=\|u_j - u_{j-1}\|_{L^s_t L^q_x(S_T)},
\]
where
\begin{equation} \label{5.32}
  \frac 1 s +\frac {\left(m+2\right)n}{ 2 q} =
  \frac{\left(m+2\right)n}{2 \delta} =\frac {m+2} 2 \left( \frac n 2 -
  \gamma\right).
\end{equation}

\smallskip

We claim that there exist a constant $\varepsilon_0>0$ and a
$\theta\in[0,1]$ such that
\begin{equation}  \label{5.34}
  2 H_0(T)^\theta \left( 2 H_0(T) + \| u_0\|_{L^\infty_t
    L^\delta_x(S_T)} \right)^{1-\theta} \le \varepsilon_0
\end{equation}
and
\begin{equation} \label{5.33}
H_j (T) \le 2  H_0(T), \qquad N_j(T) \le \frac { 1} { 2} N_{j-1}(T).
\end{equation}
Indeed, due to \eqref{5.32}, from Sobolev's embedding theorem we have
that
\begin{equation*} 
\| u (t, \cdot)\|_{L^\delta(\R^n)}    \lesssim     \| u(t, \cdot)\|_{\dot{H}^\gamma (\R^n)}.
\end{equation*}
Applying H\"older's inequality, we get that
\begin{equation*} 
 \|u_j\|_{L^{\frac {\mu_\ast (\kappa -1)} 2}(S_T)} \le \|u_j\|_{L^s_t
   L^q_x (S_T)}^\theta \|u_j\|_{L^\infty_t L^\delta_x
   (S_T)}^{1-\theta},
\end{equation*}
where $\theta = \frac 2 {n(m+2)+2} + \frac {4n(m+2)}{ \mu_\ast
  (m+2)(n-1)(q-2) + 2mq}.$ Note that $ 0 \le \theta \le 1$ for $\gamma
\ge \frac 1 {m+2}$.

\smallskip

By the same arguments as in the proof of Theorem~\ref{thm1.1}, we get
that \eqref{5.34} and \eqref{5.33} hold. Consequently, \eqref{5.2} and
\eqref{5.3} also hold. Hence, the limit $u \in
C_t^0\dot{H}_x^\gamma(S_T) \cap L^s_tL^q_x(S_T)$ of the sequence
$\{u_j\}$ is a solution of the Cauchy problem \eqref{1.1} in $S_T$.
Moreover, by Fatou's lemma and \eqref{5.33}, we have that
\begin{equation*}
  \|u\|_{C_t^0\dot{H}_x^\gamma(S_T)} + \|u\|_{L^s_t L^q_x(S_T)} \le 2
  H_0(T),
\end{equation*}
which together with \eqref{5.34} yields that $u$ satisfies estimate
\eqref{1.4}.

Further, by the same arguments as in the proof of
Theorem~\ref{thm1.1}, it follows that if both $u, \tilde u $ solve the
Cauchy problem \eqref{1.1} in $S_T$, then $u =\tilde u $ in $S_T$.
\qed


\subsection{Proof of Theorem~\ref{thm1.6}}
From the assumptions of Theorem~\ref{thm1.6}, we have
\[
  \gamma = \frac {n+1} 2 \left( \frac 1 2 -\frac 1 q \right) -\frac m
         {2 \mu_\ast (m+2)}
\]
and
\[
  -\,\frac{m}{2\mu_*\left(m+2\right)} \le \gamma < \frac 1 {m+2}
  -\frac{2(n+1)}{\mu_\ast (m+2)(n-1)} =\frac 3 {m+2} -\frac
  {n\left(2\mu_\ast -m\right)}{\mu_\ast (m+2)(n-1)}.
\]
To verify \eqref{5.2}, we set
\[
  H_j(T)=\|u_j\|_{C_t^0\dot{H}_x^\gamma(S_T)} + \|u_j\|_{L^s_t L^q_x
    (S_T)}, \quad N_j(T)= \|u_j -u_{j-1}\|_{L^s_t L^q_x(S_T)}.
\]

Let $p=q/\kappa$. Then
\[
  \frac {2n}{\left(n+1\right) p} =\frac 1 q +\frac {6\mu +m}{\mu (m+2)(n+1)}
  -\frac{n-1}{2(n+1)}.
\]
Thus we can apply Theorem~\ref{thm4.5} in case (ii) together with
H\"older's inequality to find that
\begin{multline*}
  \|u_{j+1} -u_{k+1}\|_{C_t^0\dot{H}_x^\gamma(S_T)} + \|u_{j+1}
  -u_{k+1}\|_{L^s_t L^q_x (S_T)} \\
  \lesssim \|F(u_j) -F(u_k)\|_{L^2_t L^p_x(S_T)} \lesssim \| G(u_j,
  u_k)\|_{L^\rho_t L^\sigma_x(S_T)} \| u_j -u_k\|_{L^s_t L^q_x (S_T)},
\end{multline*}
where $1/\rho =1/2 -1/s$, $1/\sigma = 1/p - 1/q=(\kappa -1)/q$.

Note that $s > (\kappa-1) \rho$ when $\gamma < \frac 1 {m+2}
-\frac{2(n+1)}{\mu_\ast (m+2)(n-1)} $. Due to condition
\eqref{1.2} and H\"older's inequality,
\begin{align*}
  \| G(u_j, u_k)\|_{L^\rho_t L^\sigma_x(S_T)} & \lesssim
  \|u_j\|_{L^{\rho(\kappa-1)}_t L^q_x (S_T)}^{\kappa -1} +
  \|u_k\|_{L^{\rho(\kappa-1)}_t L^q_x (S_T)}^{\kappa -1} \\
  & \lesssim T^{\frac 1 2 -\frac 1
    s}\bigl(\|u_j\|_{L^s_t L^q_x (S_T)}^{\kappa -1} +
  \|u_k\|_{L^s_t L^q_x (S_T)}^{\kappa -1} \bigr).
\end{align*}
As in the proof of Theorem~\ref{thm1.1}, we get that
\begin{equation} \label{5.42}
  H_j(T) \le 2 H_0(T), \quad N_j(T) \le \frac 1 2 N_{j-1}(T),
\end{equation}
and
\begin{equation} \label{5.43}
  N_0(T) \le H_0(T) T^{1/2 - \kappa /s} \le \varepsilon_0,
\end{equation}
for $\varepsilon_0>0$ small by choosing $T>0$ small. Therefore, there
is a function $u \in C_t^0\dot{H}_x^\gamma(S_T)$ $\cap \,
L^s_tL^q_x(S_T)$ such that
\[
  u_j \rightarrow u \quad \text{in $L^s_tL^q_x(S_T)$ \ as $j \rightarrow
    \infty$}
\]
and \eqref{5.2} holds. Combining Fatou's lemma and \eqref{5.42}, we
see that
\begin{equation*} 
  \|u\|_{C_t^0\dot{H}_x^\gamma(S_T)} + \|u\|_{L^s_t L^q_x (S_T)} \le 2
  H_0(T).
\end{equation*}
Together with \eqref{5.43} we get that $u$ satisfies estimate
\eqref{1.4}.

Moreover, since $2 \kappa >s$, by condition \eqref{1.2} and
H\"{o}lder's inequality, we have that, for $p=q/\kappa$,
\[
\|F(u)\|_{L^2_t L^p_x (S_T)}\lesssim \|u\|_{L^{2 \kappa}_t L^{q}_x
  (S_T)}^\kappa\lesssim T^{\frac 1 {2 } -\frac \kappa s} \|u\|_{L^s_t
  L^q_x(S_T)}^\kappa
\]
and
\begin{align*}
  \| F(u_j) -F(u)\|_{L^2_t L^p_x (S_T)} & \lesssim T^{\frac 1 {2 }
    -\frac 1 s} \bigl( \|u_j\|_{L^s_t L^q_x(S_T)}^{\kappa-1} +
  \|u\|_{L^s_t L^q_x(S_T)}^{\kappa -1} \bigr) \|u_j -u\|_{L^s_t L^q_x
    (S_T)} \\
  & \lesssim T^{\frac 1 {2 } -\frac 1 s} H_0(T)^{\kappa-1} \|u_j
  -u\|_{L^s_t L^q_x (S_T)},
\end{align*}
Therefore, $F(u) \in L^2_t L^{q/\kappa}_x (S_T)$ and $F(u_j) \to F(u)$
in $L^2_t L^{q/\kappa}_x(S_T)$ as $j \rightarrow \infty$, hence
\eqref{5.3} holds. Consequently, the limit function $u \in C_t^0
\dot{H}_x^\gamma(S_T) \cap L^s_tL^q_x(S_T)$ solves the Cauchy problem
\eqref{1.1} in $S_T$.

\smallskip

Now suppose $u, \tilde u \in C_t^0 \dot{H}_x^\gamma(S_T) \cap
L^s_tL^q_x(S_T) $ both solve the Cauchy problem \eqref{1.1} in
$S_T$. Then $v=u-\tilde u \in C_t^0 \dot{H}_x^\gamma(S_T) \cap
L^s_tL^q_x(S_T)$ is a solution of Eq.~\eqref{U1}.  Applying
Theorem~\ref{thm4.5} in case (ii) and H\"{o}lder's inequality, it
follows that
\begin{align*}
\| v \|_{L^s_t L^q_x(S_T)} & \le C\, \| G(u, \tilde u ) v \|_{L^2_t
  L^p_x(S_T)} \le C T^{\frac 1 {2 } -\frac 1 s} \bigl( \|u\|_{L^s_t
  L^q_x(S_T)}^{\kappa-1} + \|\tilde u \|_{L^s_t L^q_x(S_T)}^{\kappa
  -1} \bigr) \| v\|_{L^s_t L^q_x (S_T)} \\ & \le C T^{\frac 1 {2 }
  -\frac 1 s} H_0(T)^{\kappa-1} \| v \|_{L^s_t L^q_x (S_T)} \le
\frac 1 2\, \|v \|_{L^s_t L^q_x(S_T)}.
\end{align*}
Thus  \eqref{U2} holds and $u=\tilde u  $ in $S_T$. \qed



\end{document}